\documentclass[11pt]{article}
\usepackage{geometry}
\usepackage{graphicx}\usepackage{algorithm}\usepackage{caption}
                  
\usepackage{amssymb,amsmath,amsthm}                 
\usepackage{color,comment}

\usepackage[english]{babel}
\usepackage{url}
\usepackage{hyphenat}

\usepackage{algpseudocode}
\usepackage{multirow}
\usepackage{hyperref} 
\usepackage{tikz}
\usepackage[T1]{fontenc}
\usepackage[utf8]{inputenc}
\newtheorem{theorem}{Theorem}[section]
\newtheorem{lemma}[theorem]{Lemma}
\newtheorem{definition}[theorem]{Definition}

\title{Optimal Control for Anti-Abeta Treatment in Alzheimer's Disease using a Reaction-Diffusion Model}

 \author{
Sun Lee\thanks{Department of Mathematics, The Pennsylvania State University, State College, PA. Email: \texttt{skl5876@psu.edu}, \texttt{wxh64@psu.edu}.}
\and
Chiu-Yen Kao\thanks{Department of Mathematical Sciences, Claremont McKenna College, Claremont, CA. Email: \texttt{Chiu-Yen.Kao@ClaremontMcKenna.edu}.}
\and
Zhiyuan Li\thanks{Department of Mathematics, Ohio State University, Columbus, OH. Email: \texttt{li.9036@osu.edu}.}
\and
Tingting Dan\thanks{Department of Psychiatry, University of North Carolina at Chapel Hill, Chapel Hill, NC. Email: \texttt{tingting\_dan@med.unc.edu}.}
\and
Guorong Wu\thanks{Department of Psychiatry, Department of Computer Science, Department of Statistics and Operations Research, University of North Carolina at Chapel Hill, Chapel Hill, NC. Email: \texttt{grwu@med.unc.edu}.}
\and
Wenrui Hao\footnotemark[1]
}

\begin{document}

\maketitle
\begin{abstract}

Alzheimer’s disease (AD) is a progressive neurodegenerative disorder that severely impairs survival and quality of life. While anti-amyloid beta (A$\beta$) therapies can slow disease progression, their efficacy depends on personalized dosing that maximizes benefits and minimizes risks such as amyloid-related imaging abnormalities (ARIA). Mathematical modeling offers a powerful tool for understanding AD dynamics and optimizing treatment, yet most models focus solely on temporal behavior, overlooking spatial heterogeneity within the brain. 

In this study, we propose a spatially explicit reaction-diffusion model to describe A$\beta$ plaque dynamics. We formulate an optimal control problem to minimize plaque concentration while balancing therapeutic efficacy and treatment risk. Under reasonable assumptions, we establish well-posedness and uniqueness of the optimal solution. A Finite Element Method (FEM)-based numerical framework is developed to compute personalized treatment strategies. Our model is calibrated using longitudinal A$\beta$ positron emission tomography (PET) data from the Alzheimer’s Disease Neuroimaging Initiative (ADNI), enabling estimation of patient-specific parameters such as growth rate and effective diffusivity. 

Results show that optimized treatment strategies consistently outperform constant dosing regimens across patient groups, achieving substantial reductions in cumulative amyloid burden while minimizing side effects. This integrated, data-driven framework advances personalized, spatially informed therapeutic optimization for Alzheimer’s disease. 

\end{abstract}


\section{Introduction}
Alzheimer’s disease (AD) is one of the most prevalent and debilitating neurodegenerative disorders, affecting millions of individuals worldwide. Characterized by progressive cognitive decline, AD significantly impairs quality of life and ultimately leads to mortality, with survival ranging from 3 to 11 years after diagnosis \cite{brookmeyer2002survival,larson2004survival}. In elderly populations, AD and vascular dementia are major contributors to mortality rates, further underscoring the urgent need for effective therapeutic strategies \cite{liang2021mortality}.  

Recent advancements in anti-amyloid beta (\( A\beta \)) therapies have led to the FDA approval of several treatments aimed at slowing disease progression. {\em Aducanumab (Aduhelm)}, approved in 2021, was the first monoclonal antibody targeting amyloid-beta plaques \cite{alexander2021evaluation}. More recently, {\em Lecanemab (Leqembi)} received FDA approval in 2023 after demonstrating efficacy in reducing amyloid burden and slowing cognitive decline \cite{van2023lecanemab}. Additionally, {\em Donanemab} has shown promising results in clinical trials and is currently under FDA review \cite{mintun2021donanemab}. While these therapies represent a breakthrough in AD treatment, their effectiveness depends heavily on the {\em optimal and personalized dosing regimen}. Over-treatment can lead to adverse effects such as amyloid-related imaging abnormalities (ARIA), whereas insufficient dosing may fail to halt disease progression \cite{mayeux1999treatment,vaz2020alzheimer}. Therefore, a systematic approach to dose optimization is essential to maximize therapeutic benefits while minimizing risks.  

Mathematical modeling has emerged as a powerful tool for understanding AD progression and optimizing therapeutic interventions. Numerous mathematical frameworks have been developed to describe the causal mechanisms underlying AD, including models that incorporate temporal disease dynamics, biomarker progression, and cognitive decline \cite{bertsch2018well,hao2016mathematical,lam2022introduction,petrella2019computational,petrella2024personalized,vosoughi2020mathematical,zheng2022data}. While some models have investigated optimal treatment strategies, most focus on temporal optimization without accounting for spatial heterogeneity \cite{hao2022optimal,rabiei2025data}. However, neurodegeneration in AD is inherently spatial, as pathological changes such as amyloid plaque accumulation and tau propagation exhibit distinct regional patterns within the brain. Consequently, capturing these spatial effects is essential for developing more accurate treatment models.  

To address this gap, we propose a spatially explicit reaction-diffusion model based on the Fisher model \cite{harris2023critical}, formulated within a partial differential equation (PDE) framework. Our approach builds upon prior work on optimal treatment strategies \cite{huffaker1992optimal,lenhart1995optimal, rockne2010predicting, yousefnezhad2021optimal} by incorporating spatial dynamics and explicitly modeling the side effects associated with anti-amyloid beta (A$\beta$) treatments. Specifically, we introduce a penalty term to account for the potential adverse effects of higher treatment doses. The objective is to minimize both the amyloid-beta plaque burden and the side effects of the treatments, balancing therapeutic efficacy and patient safety.

In this study, we establish the well-posedness and uniqueness of the optimal control problem under specific assumptions. We also develop a numerical approach using the Finite Element Method (FEM) to compute the optimal treatment strategy for patients based on their amyloid-beta PET scan imaging data. Our results underscore the importance of incorporating spatial dynamics into treatment planning and offer insights into personalized therapeutic strategies for AD. By integrating PDE-based modeling with optimal control theory, this framework presents a novel approach for refining treatment protocols for emerging AD therapies.

The paper is organized as follows: In \S \ref{problem_setting}, we introduce the problem setup, preliminaries, and notation. In \S \ref{Exist}, we present the necessary conditions for optimality and establish the {existence and uniqueness} of the optimal solution for sufficiently large \( \alpha \gg 1 \), using the necessary condition computations. In \S \ref{algo}, we describe a algorithm—the Linear Combination Adjoint Method—for computing the optimal solution using the FEM. Finally, in \S \ref{Nu_Re}, we present numerical results obtained from the algorithm, using personalized PET amyloid-beta scan imaging data for different patients.

\section{Problem Formulation}\label{problem_setting}

The transport and accumulation of amyloid-beta (A$\beta$) in the brain involve several biological processes, including passive diffusion through the extracellular matrix, active clearance via the glymphatic and vascular systems, and prion-like transneuronal propagation along neural pathways. Although these mechanisms collectively shape the spatial and temporal distribution of A$\beta$, incorporating them all within a single mathematical model poses significant challenges due to limited quantitative data and the complexity of their interactions.

In this study, we approximate A$\beta$ transport using a reaction-diffusion equation with a constant diffusion coefficient, representing an effective local transport mechanism over short to intermediate timescales. This simplification, widely employed in neurodegeneration modeling, captures combined passive and active clearance effects in a tractable form \cite{Raj2012,glahn2007neuroimaging}.

To describe local plaque accumulation, we adopt a logistic growth term in Eq.~\ref{eq:model}, where the carrying capacity reflects saturation effects caused by limited available monomeric A$\beta$ and spatial constraints in the tissue microenvironment. While in reality A$\beta$ aggregation follows a nucleation-dependent polymerization process \cite{Knowles2014}, the logistic approximation serves as a practical surrogate that reproduces qualitative aggregation behavior. Future extensions of this work will incorporate more detailed nucleation-aggregation kinetics to better reflect plaque formation dynamics.

The initial distribution of A$\beta$ concentration is derived from patient-specific PET imaging data, which captures spatial heterogeneity associated with known variations in neuronal connectivity, ApoE genotype expression, and blood-brain barrier permeability \cite{thal2006}. This data-driven initialization enhances biological relevance while maintaining computational feasibility.

Formally, let \( u(\mathbf{x}, t) \) denote the concentration of A$\beta$ in the brain domain \( \Omega \subset \mathbb{R}^d \), where \( d \) is the spatial dimension. The proliferation rate is denoted by \( \rho > 0 \), and the diffusion coefficient \( D(\mathbf{x}) \in L^{\infty}(\Omega) \) satisfies \( D(\mathbf{x}) \geq \theta > 0 \) for some positive constant \( \theta \).

To model the therapeutic effect of anti-A$\beta$ treatment, we introduce a dosing function \( C(t) \), where \( C(t) \in L^{\infty}(0, T) \) represents the treatment intensity administered over the time interval \( (0, T) \). The governing reaction-diffusion model is given by \cite{rockne2010predicting}:
\begin{equation} \label{eq:model}
\begin{cases}
\displaystyle u_{t} - \nabla \cdot (D(\mathbf{x}) \nabla u) = \rho (1 - u) u - C(t) u, & \quad (\mathbf{x}, t) \in \Omega \times (0, T), \\[2ex]
\displaystyle \frac{\partial u}{\partial \mathbf{n}} = 0, & \quad (\mathbf{x}, t) \in \partial \Omega \times (0, T), \\[2ex]
u(\mathbf{x}, 0) = u_{0}(\mathbf{x}), & \quad \mathbf{x} \in \Omega.
\end{cases}
\end{equation}

Here, \( \rho \) represents the production rate of A$\beta$ plaques, and the carrying capacity is normalized to 1. The vector \( \mathbf{n} \) denotes the outward unit normal on \( \partial \Omega \), and \( u_{0}(\mathbf{x}) \) specifies the initial plaque distribution obtained from imaging data. The no-flux boundary condition reflects the anatomical constraint that A$\beta$ remains confined within the brain domain.

From a biological perspective, we aim to evaluate the effectiveness of the therapy in reducing the \( A_{\beta} \) concentration in the brain while simultaneously controlling potential side effects. To achieve this, we consider the following objective function:
\begin{equation} \label{eq:objective}
\min_{0 \leq C(t) \in L^{\infty}(0,T)} \mathcal{J}(C) = \min_{0 \leq C(t) \in L^{\infty}(0,T)}\int_{0}^{T} \left( \int_{\Omega} u_{C}(\mathbf{x}, t) \, d\mathbf{x} + \alpha C^{2}(t) \right) dt,
\end{equation}
where \( u_{C} \) is the solution of \eqref{eq:model} corresponding to the treatment function \( C(t) \). The term \( \alpha C^{2}(t) \) represents the therapy’s side effects, such as ARIA, while \( \alpha \) is a given constant and serves as a balancing coefficient that regulates the trade-off between minimizing \( A_{\beta} \) plaque burden and limiting adverse effects.

\subsection{Preliminaries and Notation}

Let \( \Omega \) be a bounded domain with a smooth boundary in \( \mathbb{R}^d \). Given a fixed time \( T \in (0, \infty) \), define \( Q_T = \Omega \times (0, T) \). For a positive integer \( k \), let \( H^k(\Omega) = W^{k,2}(\Omega) \) denote the standard Sobolev space equipped with the norm 
\begin{equation*}
    \|u\|_{H^k} = \left(\sum_{|\alpha| \leq k} \int_{\Omega} \left|D^\alpha u\right|^2 \, dx \right)^{1 / 2},
\end{equation*}
where \( \alpha = (\alpha_1, \dots, \alpha_d) \) is a multi-index with \( |\alpha| = \alpha_1 + \cdots + \alpha_d \), and 
\[
D^\alpha u = \left(\frac{\partial}{\partial x_1}\right)^{\alpha_1} \cdots \left(\frac{\partial}{\partial x_d}\right)^{\alpha_d} u.
\]

For \( p > 1 \), let \( L^p(0,T; H^k(\Omega)) \) denote the space of all functions \( u \) such that, for almost every \( t \in (0, T) \), \( u(t) = u(\cdot, t) \in H^k(\Omega) \). This space is equipped with the norm
\begin{equation*}
    \|u\|_{L^p\left(0, T; H^k(\Omega)\right)} = \left(\int_0^T \|u(t)\|_{H^k(\Omega)}^p \, dt\right)^{1 / p}.
\end{equation*}

Let \( H^1(\Omega)^* \) denote the dual space of \( H^1(\Omega) \), and \( L^p(Q_T) = L^p(0, T; L^p(\Omega)) \).

\begin{definition}
    A function \( u \in L^2(0,T; H^1(\Omega)) \) with \( u_t \in L^2(0,T; H^1(\Omega)^*) \) and \( u(\mathbf{x}, 0) = u_0(\mathbf{x}) \) is called a weak solution of \eqref{eq:model} if it satisfies weak equation:
    \begin{equation}\label{eq:weak}
        \int_{\Omega} u_t \phi \, d\mathbf{x} + \int_{\Omega} D(\mathbf{x}) \nabla u \cdot \nabla \phi \, d\mathbf{x} = \int_{\Omega} \big(\rho(1-u)u - C(t)u\big) \phi \, d\mathbf{x},
    \end{equation}
    for all test functions \( \phi \in H^1(\Omega) \) and almost every \( t \in (0, T) \).
\end{definition}

The following lemma, proved in \cite{finotti2012optimal}, establishes the existence and uniqueness of a non-negative weak solution to \eqref{eq:model}.

\begin{lemma}[Existence and Uniqueness of the weak solution] \label{thm:existence} Let \( T \in (0,\infty) \), \( C \in L^{\infty}(0,T) \) be non-negative, and \( u_0 \in L^{\infty}(\Omega) \cap H^1(\Omega) \) be non-negative. Then, for each \( C \in L^{\infty}(0,T) \), there exists a unique non-negative weak solution \( u_C \) of Eq. \eqref{eq:model}. Moreover, there exists a positive constant \( M \), depending only on \( |\Omega| \), \( T \), \( \rho \), \( \theta \), \( \|u_0\|_{L^{\infty}(\Omega)} \), \( \|C\|_{L^{\infty}(0,T)} \), and spatial dimension $d$, such that
    \begin{equation*}
        \|u_C\|_{L^{\infty}\left(Q_T\right)} \leq M.
    \end{equation*}
\end{lemma}

The next lemma, proved in \cite{yousefnezhad2021optimal}, provides a uniform bound on the solution of Eq. \eqref{eq:model}.

\begin{lemma}[Uniform Bounds] \label{thm:uniform}
    Let \( T \in (0,\infty) \), \( C \in L^{\infty}(0,T) \) be non-negative, and \( u_0 \in L^{\infty}(\Omega) \cap H^1(\Omega) \) be non-negative, and let \( u = u_C \) be the corresponding solution of  Eq. \eqref{eq:model}. Then, for each \( C \in L^{\infty}(0,T) \), there exists a positive constant \( K \), depending only on \( |\Omega| \), \( T \), \( \rho \), \( \theta \), \( \|u_0\|_{L^{\infty}(\Omega)} \), \( \|C\|_{L^{\infty}(0,T)} \), and spatial dimension $d$, such that
    \begin{equation*}  
        \|u\|_{L^{\infty}\left(0, T; L^2(\Omega)\right)} + \|u\|_{L^2\left(0, T; H^1(\Omega)\right)} + \|u_t\|_{L^2\left(0, T; H^1(\Omega)^*\right)} \leq K.
    \end{equation*}
\end{lemma}

\section{Existence and Uniqueness of the Optimal Solution}\label{Exist}
In the previous section, we established that Eq.~\eqref{eq:model} admits a unique positive solution for all non-negative initial conditions \( u_0 \in L^\infty(\Omega) \cap H^1(\Omega) \) under the given assumptions. However, the existence of a solution to Eq.~\eqref{eq:model} does not automatically guarantee the existence of an optimal control \( C \in L^\infty(0,T) \) that minimizes the objective functional in Eq.~\eqref{eq:objective}. In this section, we establish the existence of such an optimal control. 

To this end, we first examine the properties of the optimal solution and derive the necessary conditions for a minimizer~\cite{yousefnezhad2021optimal}. As a preliminary step, we differentiate the mapping \( C \in L^{\infty}(0,T) \mapsto u_C(\mathbf{x}, t) \) with respect to \( C \), which leads to the following lemma. The proof of this lemma follows directly from the arguments presented in~\cite{yousefnezhad2021optimal}.

\begin{lemma}[Sensitivity Equation] \label{thm:necessary_condition}
    Let \( C \in L^{\infty}(0,T) \) be non-negative, and let \( u = u_C \) be the corresponding solution of \eqref{eq:model}. For \( \eta \in L^{\infty}(0,T) \), the mapping \( C \to u_C(\mathbf{x}, t) \) is differentiable in the following sense: there exists a function \( \psi = \psi_{C, \eta} \in L^2(0,T; H^1(\Omega)) \) such that
    \begin{equation*}
        \psi_\epsilon \rightharpoonup \psi \quad \text{weakly in } L^2(0, T; H^1(\Omega)) \quad \text{as } \epsilon \to 0,
    \end{equation*}
    where \( \psi_\epsilon = \psi_{C, \epsilon, \eta} := \frac{u_\epsilon - u}{\epsilon} \), \( u_\epsilon := u_{C + \epsilon \eta} \), and the sensitivity \( \psi \) is the weak solution of the following {\textbf{sensitivity equation}}:
\begin{equation} \label{eq:sensitivity}
    \begin{array}{c}
        \begin{cases}
            \psi_{t} - \nabla \cdot \big( D(\mathbf{x}) \nabla \psi \big) - (\rho - 2 \rho u - C) \psi = -\eta u, 
            & \quad \text{in } \Omega \times (0, T), \\[1ex]
            \frac{\partial \psi}{\partial \mathbf{n}} = 0, 
            & \quad \text{on } \partial \Omega \times (0, T), \\[1ex]
            \psi(\mathbf{x}, 0) = 0, 
            & \quad \text{in } \Omega.
        \end{cases}
    \end{array}
\end{equation}

Thus the objective functional satisfies the expansion:
    \begin{equation} \label{eq:obj_expansion}
        \mathcal{J}(C + \epsilon \eta) = \mathcal{J}(C) + \epsilon \int_0^T \left( \int_{\Omega} \psi \, d\mathbf{x} + 2 \alpha \eta C \right) dt + o(\epsilon).
    \end{equation}
\end{lemma}
The detailed derivation is shown in Appendix~\ref{thm:necessary_condition_proof}. Next, we introduce the adjoint equation and discuss its properties in the following lemma.

\begin{lemma} \label{thm:adjoint_solution}
    Let \( C \in L^{\infty}(0,T) \) be non-negative, and let \( u = u_C \) be the corresponding solution of \eqref{eq:model}. Then, there exists \( w = w_C \in L^2(0,T; H^1(\Omega)) \) with \( w_t \in L^2(0,T; H^1(\Omega)^*) \) such that \( w \) is the weak solution of the following         {{\textbf{adjoint equation}}}:
\begin{equation} \label{eq:adjoint_equation}
    \begin{array}{c}
        \begin{cases}
            w_{t} + \nabla \cdot \big( D(\mathbf{x}) \nabla w \big) + (\rho - 2 \rho u - C) w = 1, & \text{in } \Omega \times (0, T), \\
            \frac{\partial w}{\partial \mathbf{n}} = 0, & \text{on } \partial \Omega \times (0, T), \\
            w(\mathbf{x}, T) = 0, & \text{in } \Omega.
        \end{cases}
    \end{array}
\end{equation}
    
    Moreover, there exists a constant \( K > 0 \), depending on depending only on \( |\Omega| \), \( T \), \( \rho \), \( \theta \), \( \|u_0\|_{L^{\infty}(\Omega)} \), \( \|C\|_{L^{\infty}(0,T)} \), and spatial dimension $d$, such that:
    \begin{equation*}
        \|w\|_{L^{\infty}(0, T; L^2(\Omega))} + \|w_t\|_{L^2(0, T; H^1(\Omega)^*)} + \|w\|_{L^2(0, T; H^1(\Omega))} \leq K.
    \end{equation*}
    Additionally, we can confirm \( w \leq 0 \) almost everywhere, and there exists a positive constant \( M \), depending only on \( |\Omega| \), \( T \), \( \rho \), \( \theta \), \( \|u_0\|_{L^{\infty}(\Omega)} \), \( \|C\|_{L^{\infty}(0,T)} \), and spatial dimension $d$, such that
    \begin{equation*}
        \|w\|_{L^{\infty}(Q_T)} \leq M.
    \end{equation*}
\end{lemma}
The detailed proof is shown in Appendix~\ref{thm:adjoint_solution_proof}.

\vspace{0.5em}
We now derive the necessary condition for \( C \) to be a local minimizer of the functional \( \mathcal{J}(C) \).

\begin{lemma}[Necessary Condition for Optimality]\label{Necessary_Lemma}
Let \( C \in L^{\infty}(0, T) \) be a local minimizer of \( \mathcal{J} \), and let \( \eta \in L^{\infty}(0, T) \) be arbitrary. Then:
\begin{equation} \label{eq:optimality_condition}
    \int_0^T \eta \left( 2 \alpha C + \int_{\Omega} u w \, d\mathbf{x} \right) dt = 0,
\end{equation}
where \( w \) is the weak solution of the adjoint equation \eqref{eq:adjoint_equation}.
\end{lemma}

\begin{proof}
By the definition of a local minimizer and the expansion \eqref{eq:obj_expansion}, we have:
\begin{equation*}
    0 \leq \lim_{\epsilon \to 0^+} \frac{\mathcal{J}(C + \epsilon \eta) - \mathcal{J}(C)}{\epsilon}.
\end{equation*}
From \eqref{eq:obj_expansion}, this yields:
\begin{equation*}
    0 \leq \lim_{\epsilon \to 0^+}\int_0^T \int_{\Omega} \psi_{\epsilon} \, d\mathbf{x} \, dt + \int_0^T 2 \alpha \eta C \, dt.
\end{equation*}
Using the weak formulation of the adjoint equation and integrating by parts with Eq. \eqref{eq:sensitivity}, we find:
{\small
\begin{equation}\label{eq:change}
\begin{aligned}
&\lim_{\epsilon \to 0^+}\int_0^T \int_{\Omega} \psi_{\epsilon}
=\lim_{\epsilon \to 0^+}\int_0^T \int_{\Omega} \psi_{\epsilon}\left(w_{t} + \nabla \cdot \big( D(\mathbf{x}) \nabla w \big) + (\rho - 2 \rho u - C) w\right)\\[1ex]
&=\lim_{\epsilon \to 0^+}\int_0^T \int_{\Omega} w\left(-\frac{\partial \psi_{\epsilon}}{\partial t} + \nabla \cdot \big( D(\mathbf{x}) \nabla \psi_{\epsilon} \big) + (\rho - 2 \rho u - C) \psi_{\epsilon}\right)= \int_0^T \eta \left(  \int_{\Omega} u w \, d\mathbf{x} \right) dt.
\end{aligned}
\end{equation}}

Since \( \eta \) is arbitrary we have
\begin{equation*}
    0 \leq\int_0^T \eta \left( 2 \alpha C + \int_{\Omega} u w \, d\mathbf{x} \right) dt \hbox{~and~}    0 \leq\int_0^T -\eta \left( 2 \alpha C + \int_{\Omega} u w \, d\mathbf{x} \right) dt,
\end{equation*}
and this implies \eqref{eq:optimality_condition}.
\end{proof}
Next, we show that the second-order Gâteaux derivative of the objective functional \( \mathcal{J} \) in an arbitrary direction is strictly positive. This strict convexity implies that \( \mathcal{J} \) admits at most one minimizer within the admissible set. For a detailed discussion on the relationship between convexity and uniqueness of minimizers, we refer the reader to Chapter~5 of~\cite{gelfand2000calculus}.

\begin{theorem}[Local Strict Convexity of the Objective Functional]\label{thm:loc_conv}
Let \( T \in (0, \infty) \), and suppose the control \( C(t) \) satisfies \( 0 \leq C(t) < M \) a.e. for some constant \( M > 0 \). Assume that the initial condition \( u_0 \in L^\infty(\Omega) \cap H^1(\Omega) \) is non-negative. Then, for sufficiently large \( \alpha > 0 \), the objective functional
\[
\mathcal{J}(C) = \int_{0}^{T} \left( \int_{\Omega} u_{C}(\mathbf{x}, t) \, d\mathbf{x} + \alpha C^2(t) \right) \, dt
\]
is strictly convex with respect to \( C \). Moreover, there exists a constant \( \alpha_0 > 0 \), depending only on \( |\Omega| \), \( T \), \( \rho \), \( \theta \), \( \|u_0\|_{L^{\infty}(\Omega)} \), \( \|C\|_{L^{\infty}(0,T)} \), and spatial dimension $d$, such that for all \( \alpha > \alpha_0 \), second-order Gˆateaux derivative of $\mathcal{J}$ strictly positive.

\end{theorem}

\begin{proof}
To verify the local strict convexity of \( \mathcal{J}(C) \), we compute its second-order Gâteaux derivative with respect to the control \( C \). Our goal is to show that the second variation is strictly positive when \( \alpha \) is sufficiently large. First, we examine the derivative of the mapping \( C \mapsto u_C(\mathbf{x}, t) \), followed by the derivation of its second-order variation. The governing state equation is given by:

\begin{equation} \label{eq:state}
\begin{cases}
    u_t - \nabla \cdot \big( D(\mathbf{x}) \nabla u \big) - (\rho - \rho u - C) u = 0, & (\mathbf{x}, t) \in \Omega \times (0, T), \\
    \frac{\partial u}{\partial \mathbf{n}} = 0, & (\mathbf{x}, t) \in \partial \Omega \times (0, T), \\
    u(\mathbf{x}, 0) = u_0(\mathbf{x}), & \mathbf{x} \in \Omega.
\end{cases}
\end{equation}

The derivatives are formulated as follows:

\begin{itemize}
    \item \textbf{First-order equation with direction \( \eta_1 \):}
        \begin{equation} \label{eq:sensitiv_1}
    \begin{cases}
        \psi_{1,t} - \nabla \cdot \big( D(\mathbf{x}) \nabla \psi_1 \big) - (\rho - 2 \rho u - C) \psi_1 = -\eta_1 u, & \text{in } \Omega \times (0, T), \\
        \frac{\partial \psi_1}{\partial \mathbf{n}} = 0, & \text{on } \partial \Omega \times (0, T), \\
        \psi_1(\mathbf{x}, 0) = 0, & \mathbf{x} \in \Omega,
    \end{cases}
    \end{equation}
where \( \psi_1 = \psi_{C, \eta_1} := \lim_{\epsilon \to 0^+}\frac{u_\epsilon - u}{\epsilon} \), with \( u_\epsilon := u_{C + \epsilon \eta_1} \).

\item \textbf{Second-order equation with direction \( \eta_2 \):}
   \begin{equation} \label{eq:sensitiv_2}
    \begin{cases}
        \psi_{2,t} - \nabla \cdot \big( D(\mathbf{x}) \nabla \psi_2 \big) - (\rho - 2 \rho u - C) \psi_2 = - (\eta_1\tilde{\psi_1}  + \eta_2\psi_1 + 2 \rho \psi_1 \tilde{\psi_1}), & \text{in } \Omega \times (0, T), \\
        \frac{\partial \psi_2}{\partial \mathbf{n}} = 0, & \text{on } \partial \Omega \times (0, T), \\
        \psi_2(\mathbf{x}, 0) = 0, & \mathbf{x} \in \Omega.
    \end{cases}
    \end{equation}
where \( \psi_2 = \psi_{C, \eta_1, \eta_2} := \lim_{\delta \to 0^+}\frac{\psi_{1,\delta} - \psi_1}{\delta} \), with \( \psi_{1,\delta} := \psi_{C + \delta \eta_2, \eta_1} \). Here, \( \tilde{\psi_1} \) represents the derivative of \( u \) in the \( \eta_2 \) direction:
\begin{equation}
    \tilde{\psi_1} = \lim_{\delta \to 0^+}\frac{u_{C+\delta \eta_2}-u_C}{\delta}.
\end{equation}
\end{itemize}

We can derive second-order equation Eq.~\eqref{eq:sensitiv_2} using similar knowledge as in Lemma \ref{eq:sensitivity}. To determine the necessary conditions for optimality, we compute the first Gâteaux derivative of \( \mathcal{J} \) using Eq.~\eqref{eq:obj_expansion}, \eqref{eq:change} as:
\begin{equation*}
    \lim_{\epsilon \to 0^+} \frac{\mathcal{J}(C + \epsilon \eta_1) - \mathcal{J}(C)}{\epsilon} = \int_{0}^{T} \int_{\Omega} \psi_1 \, d\mathbf{x} dt + \int_{0}^{T} 2 \alpha \eta_1 C \, dt.
\end{equation*}
The second-order Gâteaux derivative is given by:
\begin{align*}
    &\lim_{\delta\to 0^+} \frac{1}{\delta}\left(\lim_{\epsilon \to 0^+} \frac{\mathcal{J}(C + \epsilon \eta_1+\delta \eta_2) - \mathcal{J}(C+\delta \eta_2)}{\epsilon} - \lim_{\epsilon \to 0^+} \frac{\mathcal{J}(C + \epsilon \eta_1) - \mathcal{J}(C)}{\epsilon}\right) \\
    &= \int_{0}^{T} \int_{\Omega} \psi_2 \, d\mathbf{x} dt + \int_{0}^{T} 2 \alpha \eta_1 \eta_2 \, dt= \int_{0}^{T} \int_{\Omega} (\eta_1 \tilde{\psi_1} + \eta_2 \psi_1 + 2 \rho \psi_1 \tilde{\psi_1}) w \, d\mathbf{x} dt + \int_{0}^{T} 2 \alpha \eta_1 \eta_2 \, dt.
\end{align*}

The last term was derived using a similar technique in Eq. \eqref{eq:change}. In particular, when \( \eta_1 = \eta_2 \), we have:
\begin{equation}\label{second_derivative}
    \int_{0}^{T} \int_{\Omega} (2 \eta_1 + 2 \rho \psi_1) \psi_1 w \, d\mathbf{x} dt + \int_{0}^{T} 2 \alpha \eta_1^2 \, dt.
\end{equation}
Our goal is to demonstrate that Eq.~\eqref{second_derivative} is positive when \( \alpha \gg 1 \). Using Cauchy's inequality with Lemma \ref{thm:adjoint_solution} we get

\begin{equation}\label{eq:upper_1term}
    \int_0^T \int_\Omega |\eta_1 \psi_1 w| \, d\mathbf{x} dt \leq \frac{1}{2} |\Omega| \int_0^T \eta_1^2 \, dt + \frac{1}{2} \|w\|^2_{L^{\infty}(Q_T)} \|\psi_1\|_{L^2(Q_T)}^2,
\end{equation}
\begin{equation}\label{eq:upper_2term}
    \int_0^T \int_\Omega |\psi_1 \psi_1w| \, dx dt
\leq \|w\|_{L^{\infty}(Q_T)}   \|\psi_1\|_{L^2(Q_T)}^2.
\end{equation}
We know from Lemma~\ref{thm:adjoint_solution} that $\|w\|_{L^{\infty}(Q_T)}$ is bounded. To properly estimate Eq.~\eqref{second_derivative} value, we need to compute $ \|\psi_1\|_{L^2(Q_T)}^2$. We begin this computation by considering Eq.~\eqref{eq:sensitiv_1}. Given that \( 2 \rho u + C \geq 0 \), multiplying the Eq.~\eqref{eq:sensitiv_1} by \( \psi_1 \) and integrating over \( \Omega \times (0, t) \), we obtain:
\begin{equation} \label{eq:psi_bound_1}
    \frac{1}{2} \|\psi_1(\mathbf{x}, t)\|_{L^2(\Omega)}^2 + \int_0^t \int_{\Omega} D(\mathbf{x}) |\nabla \psi_1(\mathbf{x}, s)|^2 \, d\mathbf{x} ds 
    \leq \rho \int_0^t \|\psi_1\|_{L^2(\Omega)}^2 \, ds + \|u\|_{L^{\infty}(Q_T)} \int_0^t \int_{\Omega} |\eta_1(s) \psi_1(\mathbf{x}, s)| \, d\mathbf{x} ds.
\end{equation}

Using Cauchy's inequality, the second term satisfies:
\begin{equation} \label{eq:psi_bound_2}
    \int_0^t \int_{\Omega} |\eta_1(s) \psi_1(\mathbf{x}, s)| \, d\mathbf{x} ds \leq \frac{1}{2} |\Omega| \int_0^t \eta_1^2(s) \, ds + \frac{1}{2} \int_0^t \|\psi_1(\mathbf{x}, s)\|_{L^2(\Omega)}^2 \, ds.
\end{equation}

Substituting \eqref{eq:psi_bound_2} into \eqref{eq:psi_bound_1}, we obtain:
\begin{align}\label{eq:clean1}
    \frac{1}{2} \|\psi_1(\mathbf{x}, t)\|_{L^2(\Omega)}^2 + \int_0^t \int_{\Omega} D(\mathbf{x}) |\nabla \psi_1(\mathbf{x}, s)|^2 \, d\mathbf{x} ds 
    &\leq \rho \int_0^t \|\psi_1(\mathbf{x}, s)\|_{L^2(\Omega)}^2 \, ds + \frac{|\Omega|}{2} \|\eta_1\|^2_{L^2(0, t)} \|u\|_{L^{\infty}(Q_T)} \notag \\
    &\quad + \frac{1}{2} \|u\|_{L^{\infty}(Q_T)} \int_0^t \|\psi_1(\mathbf{x}, s)\|_{L^2(\Omega)}^2 \, ds.
\end{align}

Simplifying, we get:
\begin{equation}\label{eq:clean2}
    \|\psi_1(\mathbf{x}, t)\|_{L^2(\Omega)}^2 \leq B_1 \int_0^t \|\psi_1(\mathbf{x}, s)\|_{L^2(\Omega)}^2 \, ds + \|\eta_1\|^2_{L^2(0,t)} B_2,
\end{equation}
where \( B_1 = 2\rho + \|u\|_{L^{\infty}(Q_T)} \) and \( B_2 = |\Omega| \|u\|_{L^{\infty}(Q_T)} \).

Applying Grönwall's inequality \cite{evans2022partial}, we derive the bound:
\begin{equation} \label{eq:gronwall_bound}
    \|\psi_1\|_{L^2(Q_T)}^2 \leq \|\eta_1\|^2_{L^2(0,T)} B_2 \left( 1 + T e^{B_1T} \right).
\end{equation}

{Combining Eqs.~\eqref{eq:upper_1term},~\eqref{eq:upper_2term} and \eqref{eq:gronwall_bound}, we have}
\begin{equation}
\begin{aligned}
&\int_{0}^{T} \int_{\Omega} (2 \eta_1 + 2 \rho \psi_1) \psi_1 w \, d\mathbf{x} dt + \int_{0}^{T} 2 \alpha \eta_1^2 \, dt
\\
&\geq    \int_{0}^{T} 2 \alpha \eta_1^2 \, dt- 2 \int_0^T \int_\Omega |\eta_1 \psi_1 w| \, d\mathbf{x} dt-2\int_0^T \int_\Omega \rho |\psi_1 \psi_1w| \, dx dt
    \\
    &\geq \left(2\alpha -{|\Omega|}\right)\|\eta_1\|^2_{L^2(0,T)}-\left(\|w\|_{L^{\infty}(Q_T)}+2\rho   \right)\|w\|_{L^{\infty}(Q_T)}\|\psi_1\|_{L^2(Q_T)}^2
    \\
    &\geq  \left(2\alpha -{|\Omega|}-\left(\|w\|_{L^{\infty}(Q_T)}+2\rho   \right)\|w\|_{L^{\infty}(Q_T)} B_2 \left( 1 + T e^{B_1T} \right)\right)\|\eta_1\|^2_{L^2(0,T)}
\end{aligned}
\end{equation}
For sufficiently large 
\begin{equation}\label{inequality_convexity}
    \alpha> \left({|\Omega|}+\left(\|w\|_{L^{\infty}(Q_T)}+2\rho   \right)\|w\|_{L^{\infty}(Q_T)} B_2 \left( 1 + T e^{B_1T} \right)\right)/2,
\end{equation} this inequality ensures the positivity of Eq.~\eqref{second_derivative}.

Finally, from Lemma \ref{thm:existence} and \ref{thm:adjoint_solution} we observe that both \( \|u\|_{L^\infty(Q_T)} \) and \( \|w\|_{L^\infty(Q_T)} \) upper bound depend continuously on \( |\Omega| \), \( T \), \( \rho \), \( \theta \), \( \|u_0\|_{L^{\infty}(\Omega)} \), \( \|C\|_{L^\infty(0,T)} \), and spatial dimension $d$. Since the control is assumed to satisfy \( \|C\|_{L^\infty(0,T)} < M \), there exists a constant \( K_M > 0 \) such that
\[
\|u\|_{L^\infty(Q_T)},\, \|w\|_{L^\infty(Q_T)} \leq K_M.
\]
Therefore, the lower bound in Eq.~\eqref{inequality_convexity} depends only on \( |\Omega| \), \( T \), \( \rho \), \( \|u_0\|_{L^{\infty}(\Omega)} \), \( M \), and the spatial dimension \( d \). Consequently, the second variation of \( \mathcal{J} \) remains strictly positive for sufficiently large \( \alpha \). This implies that \( \mathcal{J} \) is strictly convex in a neighborhood of any admissible control with \( \|C\|_{L^\infty(0,T)} < M \). Hence, the objective functional is locally strictly convex when $\alpha > \alpha_0$
\[
\alpha_0=\sup_{\|C\|_{L^\infty(0,T)<M}}
\left({|\Omega|}+\left(\|w\|_{L^{\infty}(Q_T)}+2\rho   \right)\|w\|_{L^{\infty}(Q_T)} B_2 \left( 1 + T e^{B_1T} \right)\right)/2.
\]
\end{proof}

Finally, we establish the existence and uniqueness of the optimal solution for sufficiently large \( \alpha \), under the assumption that \( \|C\|_{L^\infty(0, T)} \) is bounded. We construct a fixed-point iteration based on the necessary optimality condition derived in \S~\ref{Exist} and demonstrate that it defines a contraction mapping in \( L^\infty(0,T) \) for sufficiently large \( \alpha \). By the Banach fixed-point theorem, this guarantees the existence of a unique fixed point satisfying the necessary condition for optimality. Combining this result with the strict convexity of \( \mathcal{J} \), we conclude that this fixed point corresponds to the unique optimal control. Following this reasoning, we summarize the result in the following theorem.

\begin{theorem}[Existence and Uniqueness of the Optimal Control]\label{thm:contraction}
Let \( T \in (0, \infty) \), and suppose the control \( C(t) \) satisfies \( 0 \leq C(t) < M \) almost everywhere on \( (0, T) \), for some constant \( M > 0 \). Assume that the initial condition \( u_0 \in L^\infty(\Omega) \cap H^1(\Omega) \) is non-negative. Then, for sufficiently large \( \alpha > 0 \), the objective functional
\[
\mathcal{J}(C) = \int_{0}^{T} \left( \int_{\Omega} u_{C}(\mathbf{x}, t) \, d\mathbf{x} + \alpha C^2(t) \right) \, dt
\]
admits a unique local minimizer \( C^* \in L^\infty(0, T) \) within the admissible set 
\[
\mathcal{A} = \{\, C \in L^\infty(0, T) : \|C\|_{L^\infty(0,T)} < M \,\}.
\]
\end{theorem}

\begin{proof}
First, we define a contraction mapping \( \Phi: L^\infty(0, T) \to L^\infty(0, T) \) as
\[
\Phi(C)(t) := \beta C(t) + (1 - \beta) \left( -\frac{1}{2\alpha} \int_\Omega u_C(\mathbf{x}, t)\, w_C(\mathbf{x}, t) \, d\mathbf{x} \right),
\]
for some \( \beta \in [0, 1) \).
For \( C_1, C_2 \in L^\infty(0, T) \), we have
\[
\Phi(C_2)(t) - \Phi(C_1)(t) = \beta (C_2(t) - C_1(t)) + \frac{1 - \beta}{2\alpha} \Big[ \int_\Omega u_{C_1} w_{C_1} \, d\mathbf{x} - \int_\Omega u_{C_2} w_{C_2} \, d\mathbf{x} \Big],
\]
with
\[
\int_\Omega u_{C_1} w_{C_1} \, d\mathbf{x} - \int_\Omega u_{C_2} w_{C_2} \, d\mathbf{x} = \int_\Omega \delta u w_{C_2} \, d\mathbf{x} + \int_\Omega u_{C_1} \delta w \, d\mathbf{x},
\]
where \( \delta u := u_{C_2} - u_{C_1} \) and \( \delta w := w_{C_2} - w_{C_1} \).
The sensitivity equations for \( \delta u \) and \( \delta w \) satisfy
\[
\begin{cases}
\delta u_t - \nabla \cdot \left( D(\mathbf{x}) \nabla \delta u \right) - (\rho - \rho u_1 - \rho u_2 - C_2) \delta u = (C_2 - C_1) u_1, & \text{in } \Omega \times (0, T), \\
\partial_{\mathbf{n}} \delta u = 0, & \text{on } \partial \Omega \times (0, T), \\
\delta u(\mathbf{x}, 0) = 0, & \text{in } \Omega,
\end{cases}
\]
and
\[
\begin{cases}
\delta w_t + \nabla \cdot \left( D(\mathbf{x}) \nabla \delta w \right) + (\rho - 2\rho u_2 - C_2) \delta w = \left( 2\rho(u_2 - u_1) + C_2 - C_1 \right) w_1, & \text{in } \Omega \times (0, T), \\
\partial_{\mathbf{n}} \delta w = 0, & \text{on } \partial \Omega \times (0, T), \\
\delta w(\mathbf{x}, T) = 0, & \text{in } \Omega.
\end{cases}
\]

Applying standard energy estimates and Grönwall's inequality (Lemma~\ref{thm:adjoint_solution}, Eq.~\eqref{eq:gronwall_bound}) gives
\[
\| \delta u \|_{L^\infty(0, T; L^2(\Omega))} \leq K_1 \| C_2 - C_1 \|_{L^\infty(0, T)} \| u_{C_1} \|_{L^\infty(0, T; L^2(\Omega))},
\]
\[
\| \delta w \|_{L^\infty(0, T; L^2(\Omega))} \leq K_2 \Big( \| C_2 - C_1 \|_{L^\infty(0, T)} + 2\rho \| \delta u \|_{L^\infty(0, T; L^2(\Omega))} \Big) \| w_{C_1} \|_{L^\infty(0, T; L^2(\Omega))},
\]
where \( K_1, K_2 > 0 \) depend only on \( |\Omega|, T, \rho, \theta, M \), and the spatial dimension \( d \).

Using the a priori bounds \( \| u_C \|_{L^\infty(0, T; L^2(\Omega))} \leq U_M \) and \( \| w_C \|_{L^\infty(0, T; L^2(\Omega))} \leq W_M \) for \( \|C\|_{L^\infty(0,T)} < M \) (Lemmas~\ref{thm:existence} and~\ref{thm:adjoint_solution}), we obtain
\begin{align*}
\left| \Phi(C_2)(t) - \Phi(C_1)(t) \right| 
&\leq \beta \| C_2 - C_1 \|_{L^\infty(0, T)} \\
&\quad + \frac{1 - \beta}{2\alpha} \Big( \| \delta u(t) \|_{L^2(\Omega)} \| w_{C_2}(t) \|_{L^2(\Omega)} + \| u_{C_1}(t) \|_{L^2(\Omega)} \| \delta w(t) \|_{L^2(\Omega)} \Big) \\
&\leq \left( \beta + \frac{(1 - \beta) K}{\alpha} \right) \| C_2 - C_1 \|_{L^\infty(0, T)},
\end{align*}
with \( K := (K_1 + K_2(1 + 2\rho K_1 U_M)) U_M W_M \).

Hence, for \( \alpha > \alpha_1 := \max(\alpha_0, K) \), where \( \alpha_0 \) ensures strict convexity of \( \mathcal{J} \) (Theorem~\ref{thm:loc_conv}), we have \( \gamma := \beta + \frac{(1-\beta) K}{\alpha} < 1 \), so \( \Phi \) is a contraction.

Moreover, \( \Phi \) maps the admissible set 
\[
\mathcal{U}_M := \{ C \in L^\infty(0,T) : \|C\|_{L^\infty(0,T)} < M \}
\]
into itself. By the Banach fixed-point theorem, \( \Phi \) admits a unique fixed point \( C^* \in \mathcal{U}_M \) with \( \Phi(C^*) = C^* \).

By construction, \( C^* \) satisfies the first-order necessary optimality condition
\[
\mathcal{J}'(C^*)[\eta] = \frac{d}{d\epsilon} \mathcal{J}(C^* + \epsilon \eta)\Big|_{\epsilon = 0} = 0 \quad \forall \eta \in L^\infty(0, T),
\]
and, by Theorem~\ref{thm:loc_conv}, \( \mathcal{J} \) is strictly convex over \( \mathcal{U}_M \). Therefore, \( C^* \) is the unique minimizer of \( \mathcal{J} \) in \( \mathcal{U}_M \), i.e., the unique optimal control.
\end{proof}

\section{Numerical Methods}\label{algo}
In this section, we present a numerical method to solve the optimal control problem based on the adjoint formulation derived in \S~\ref{Exist}. The optimality condition shown in Eq.~\eqref{eq:optimality_condition} provides an updated version of the optimal solution \( C \). To solve this problem, we first employ the FEM to discretize both the state equation (\ref{eq:model}) for \( u \) and the adjoint equation (\ref{eq:adjoint_equation}) for \( w \) in the spatial domain.

For simplicity, we will use the \( u \) as an example to illustrate the discretization process. We begin by defining a triangulation of the domain \( \Omega \), denoted by \( \mathcal{T}_h = \{T\}_{k=1, \dots, N_e} \). As usual, we define the mesh size \( h \) as:
\[
h = \max_{T \in \mathcal{T}_h} \text{diam}(T),
\]
where \( \text{diam}(T) \) is the diameter of element \( T \in \mathcal{T}_h \).

Let \( V_h \) be the subspace of \( L^\infty(\Omega) \cap H^1(\Omega) \) composed of piecewise globally continuous polynomials of degree \( r \geq 1 \). We define \( V_h \) as:
\[
V_h = \text{span} \{ \phi_h^k : k = 1, \dots, N_h \},
\]
where \( \phi_h^k \) is the basis of \( V_h \), and \( N_h \) is the dimension of \( V_h \). For any \( u_h \in V_h \), there exists a unique vector \( \vec{u}_h = (u_h^1, \dots, u_h^{N_h})^T \in \mathbb{R}^{N_h} \) such that:
\[u_h = \sum_{k=1}^{N_h} u_h^k \phi_h^k.\]
Discretizing the spatial domain leads to the following weak formulation of Eq.~\eqref{eq:weak}:  
\begin{equation}
     \int_{\Omega} \frac{\partial u_h}{\partial t} \phi_h \, d\mathbf{x} 
     = -\int_{\Omega} D(\mathbf{x}) \nabla u_h \cdot \nabla \phi_h \, d\mathbf{x}
       + \int_{\Omega} \left(\rho (1 - u_h)u_h - Cu_h\right) \phi_h \, d\mathbf{x}, \quad \forall \phi_h \in V_h.
\end{equation}

For temporal discretization, we treat the spatially discretized formulation as a system of ordinary differential equations (ODEs) in time:  
\begin{equation}
M\frac{dU}{dt} = F(U) \hbox{~and~ }M_{i,j} = \int_{\Omega} \phi^i_h(x) \phi^j_h(x) \, d\mathbf{x}.
\end{equation}
where  \( M = [M_{i,j}] \) is the mass matrix.
Here, \( U \) represents the vector of nodal values of \( u_h \), and \( F(U) \) denotes the spatially discretized nonlinear operator. 

\subsection{Linear Combination Adjoint Method}

We update the control \(C\) using the optimality condition, employing a linear combination of the previous control and an intermediate control computed via the adjoint method.
\begin{algorithm}[H]
\caption{Linear Combination Adjoint Method} \label{alg:linear_combination}
\hspace*{\algorithmicindent} \textbf{Input:} Parameters \( \rho, D(\mathbf{x}), \alpha, \text{TOL} > 0 \), and an initial control \( C_0 \).
\begin{algorithmic}[1]
    \State Set \( i = 0 \).
    \State Compute the state \( u_i = u_{C_i} \) by solving the state equation \eqref{eq:model} using FEM.
    \State Set \( \rho - 2 \rho u_i - C_i \) as the coefficient in \eqref{eq:adjoint_equation} and compute the adjoint state \( w_i = w_{C_i} \) using FEM.
    \State Compute the intermediate control:
    \begin{equation*}
        \tilde{C} = - \frac{1}{2\alpha} \int_{\Omega} u_i w_i \, d\mathbf{x}.
    \end{equation*}
    \State Update the control:
    \begin{equation*}
        C_{i+1} = \beta C_i + (1 - \beta) \tilde{C}, \quad \beta \in [0, 1).
    \end{equation*}
    \State If \( \|C_{i+1} - C_i\| < \text{TOL} \), stop; otherwise, set \( i \leftarrow i + 1 \) and return to Step 2.
\end{algorithmic}
\end{algorithm}

\subsection{Convergence of the Linear Combination Adjoint Method}

We now justify the convergence of Algorithm~\ref{alg:linear_combination}, based on the first-order necessary optimality condition (Lemma~\ref{Necessary_Lemma}). At the optimal control \( C^\ast \), the Gâteaux derivative of the objective functional \( \mathcal{J} \) vanishes:
\[
2\alpha C^\ast + \int_\Omega u_{C^\ast} w_{C^\ast} \, d\mathbf{x} = 0,
\]
which gives the pointwise relation
\[
C^\ast = - \frac{1}{2\alpha} \int_\Omega u_{C^\ast} w_{C^\ast} \, d\mathbf{x}.
\]

Motivated by this identity, we define the intermediate control at iteration \( i \) as
\[
\tilde{C} = - \frac{1}{2\alpha} \int_\Omega u_i w_i \, d\mathbf{x},
\]
where \( u_i = u_{C_i} \) and \( w_i = w_{C_i} \) are the state and adjoint state corresponding to \( C_i \). The updated control is obtained via the linear combination
\[
C_{i+1} = \beta C_i + (1 - \beta) \tilde{C}, \quad \beta \in [0, 1),
\]
which ensures smooth progression toward the optimal solution without overshooting.

Since, when \( \|C\|_{L^\infty(0,T)} < M \), \( \mathcal{J} \) is strictly convex (Theorem~\ref{thm:loc_conv}) and the mapping \( C \mapsto -\frac{1}{2\alpha} \int_\Omega u_C w_C \, d\mathbf{x} \) is continuous in \( L^2(\Omega) \), the fixed-point iteration
\[
C_{i+1} = \Phi(C_i) := \beta C_i + (1 - \beta) \left(-\frac{1}{2\alpha} \int_\Omega u_i w_i \, d\mathbf{x} \right)
\]
defines a contraction mapping (Theorem~\ref{thm:contraction}). Consequently, the sequence \( \{ C_i \} \) converges strongly in \( L^2(\Omega) \) to \( C^\ast \), provided the initial guess satisfies \( \|C\|_{L^\infty(0,T)} < M \) and \( \beta \) is chosen appropriately.
Moreover, if \( C_i \geq 0 \) and \( u_i \geq 0 \), \( w_i \leq 0 \) a.e., then \( \tilde{C} \geq 0 \). The convex combination \( C_{i+1} = \beta C_i + (1 - \beta)\tilde{C} \) with \( \beta \in [0,1) \) ensures \( C_{i+1} \geq 0 \), preserving admissibility at each iteration.

{\bf Remark:} When \( \beta \) is close to 1, the fixed-point iteration can be interpreted as a gradient descent scheme for minimizing \( \mathcal{J} \):
\[
C_{i+1} = C_i - \gamma \left( 2\alpha C_i + \int_\Omega u_i w_i \, d\mathbf{x} \right), \quad \gamma = \frac{1-\beta}{2\alpha} > 0,
\]
with \( u_i = u_{C_i} \) and \( w_i = w_{C_i} \) corresponding to \( C_i \). Using the first-order expansion of \( \mathcal{J} \) (Eq.~\eqref{eq:obj_expansion} and Lemma~\ref{Necessary_Lemma}),
\[
\mathcal{J}(C_i + \gamma \eta_i) = \mathcal{J}(C_i) - \gamma \int_0^T \left\| 2\alpha C_i + \int_\Omega u_i w_i \, d\mathbf{x} \right\|^2_{L^2(\Omega)} dt + o(\gamma),
\]
where \( \eta_i = - \left( 2\alpha C_i + \int_\Omega u_i w_i \, d\mathbf{x} \right) \) is the negative gradient direction.  

Due to strict convexity (Theorem~\ref{thm:loc_conv}), if \( C_i \) is not optimal, the integrand is strictly positive, implying
\[
\mathcal{J}(C_{i+1}) < \mathcal{J}(C_i),
\]
so the sequence \( \{ \mathcal{J}(C_i) \} \) decreases strictly at each iteration.

\section{Numerical Results}\label{Nu_Re}
In this section, we present numerical experiments for both 1D and 2D cases using specified parameters to verify whether the recovered results are consistent with the expected ground truth. The 2D domain is taken as a slice from the 3D PET scan data, which allows us to assess the method’s performance in a more complex geometric setting. Moreover, for the 2D experiments, we also set the initial condition from the PET amyloid data, while other parameters are given. Finally, we apply our method to the full 3D case with real data to calibrate the model parameters and implement the optimal control.

\subsection{1D case}\label{1_D_S}
\begin{figure}
    \centering
    \includegraphics[width=1\linewidth]{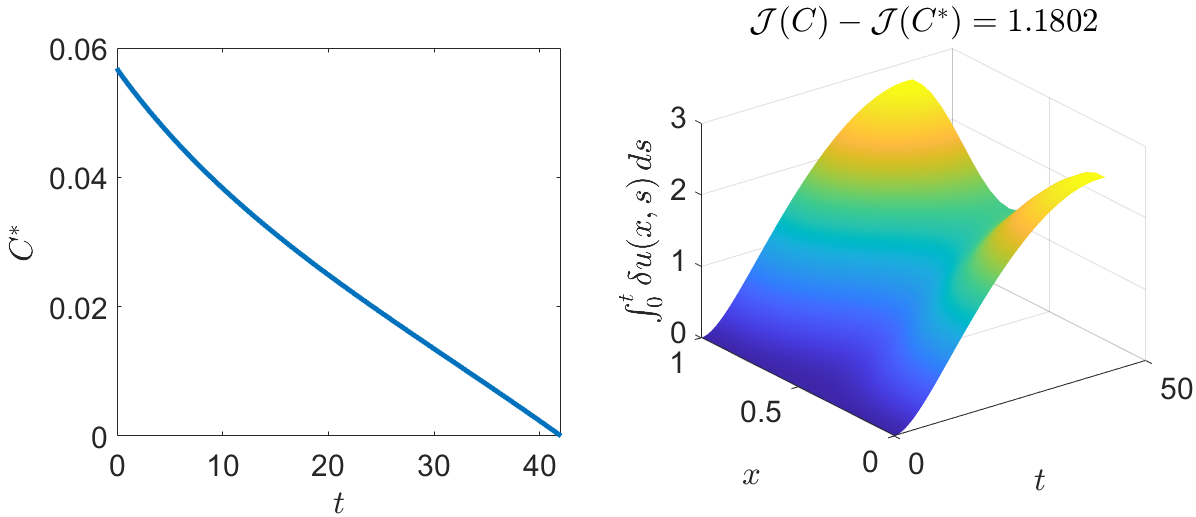}
    \caption{The difference defined as 
$\delta u(\mathbf{x}, t) = u(\mathbf{x}, t) - u^*(\mathbf{x}, t)$, 
where $u$ denotes the solution under the constant treatment $C$, 
and $u^*$ under the optimal treatment $C^*$. We set \( T = 42 \) and initialized \( u_0 \) using the PET images with \( \alpha = 100 \), \( \rho = 0.012 \), and \( D = 0.002 \) for Eq.~\eqref{eq:model} and Eq.~\eqref{eq:objective}.
}
    \label{1D}
\end{figure}

We set the initial condition as  
\[
u_0(\mathbf{x}) = \frac{\cos(2\pi x) + 1}{2}
\]
on the domain \( \Omega = [0,1] \).  
We compare the optimal treatment strategy, \( C^*(t) \), with a constant treatment, defined as  
\begin{equation}\label{c_avf}
    C = \frac{\int_0^T C^*(t) dt}{T}.
\end{equation}
We use the linear combination adjoint method and choose \( \beta = 0.5 \) in Algorithm \ref{alg:linear_combination}. The results are shown in Fig.~\ref{1D} for $\alpha = 100$ with $D = 0.002$ and $\rho = 0.012$ and \( T = 42 \) from Eq.~\eqref{eq:model}.

As illustrated in Fig.~\ref{1D}, the objective function satisfies \(\mathcal{J}(C^*(t)) < \mathcal{J}(C),\) indicating that the optimal control outperforms the constant treatment. Moreover, the necessary condition in Lemma~\ref{Necessary_Lemma} for optimality is nearly zero ($7.85\times10^{-7}$) for \(C^*(t)\), whereas it deviates significantly from zero for the constant treatment, confirming the validity of the computed optimal control. Figure~\ref{1D} also presents the temporal evolution of the cumulative integral
\[
\int_0^t \big(u(\mathbf{x},s) - u^*(\mathbf{x},s)\big)\, ds,
\]
where \(u\) denotes the solution under the constant treatment \(C\) and \(u^*\) corresponds to the solution under the optimal treatment \(C^*\). The cumulative difference steadily accumulates over time and remains strictly positive throughout the treatment period $(t>0)$, indicating a sustained advantage of the optimal strategy. By explicitly accounting for side effects within the optimization process, the optimal control achieves a superior overall balance between amyloid plaque clearance and treatment safety.

\subsection{2D Case with \texorpdfstring{$\beta$}{Beta}-Amyloid PET Imaging Data}\label{2_D_S}

Next, we apply real patient data from positron emission tomography (PET) scans of the brain, obtained from the Alzheimer's Disease Neuroimaging Initiative (ADNI) dataset \cite{mueller2005alzheimer}  (https:// adni.loni.usc.edu). The ADNI dataset is a comprehensive and widely used collection of longitudinal clinical, imaging, genetic, and other biomarker data. To illustrate the concept of optical control treatment aimed at clearing \(\beta\)-amyloid plaques in the brain, we use data from one subject in each of the following five diagnostic groups: Cognitively Normal (CN), Significant Memory Complaint (SMC) but clinically normal, Early Mild Cognitive Impairment (EMCI), Late Mild Cognitive Impairment (LMCI), and Alzheimer's Disease (AD).

\begin{figure}
    \centering
    \includegraphics[width=1\linewidth]{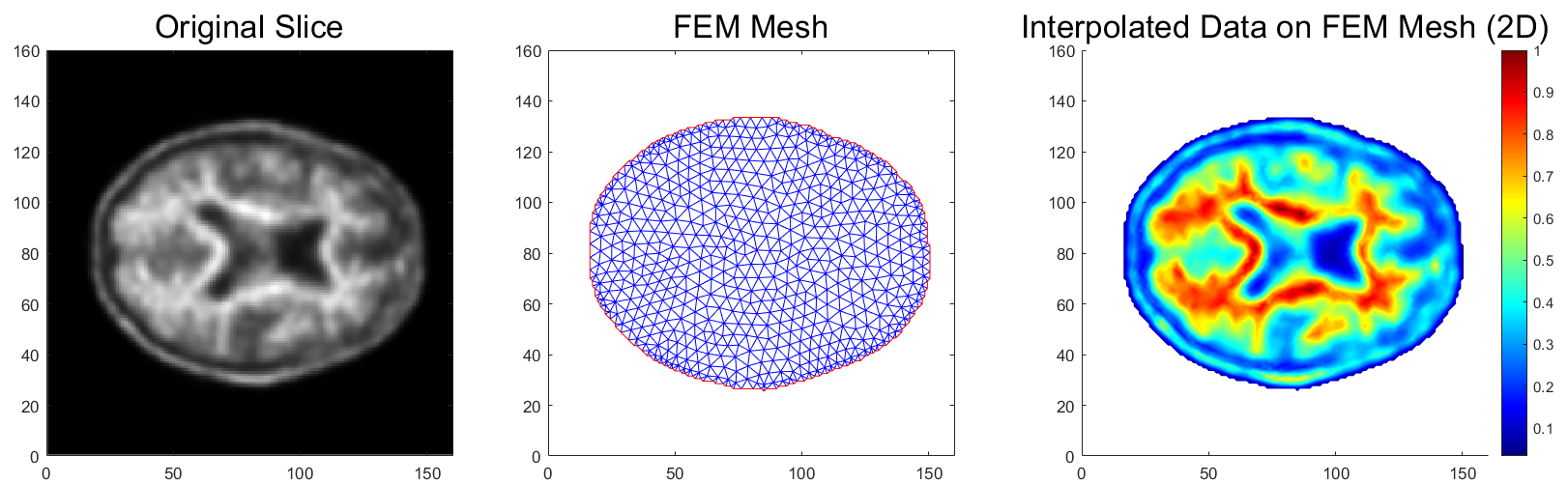}
    \caption{{\bf Left}: The original PET scan of the brain for the AD patient. {\bf Middle}: The mesh was generated based on the PET scan data. {\bf Right}: The initial condition \(u_0(x,y)\) generated from the PET scan data using the formula in Eq. \ref{IC}.}
    \label{mesh}
\end{figure}

The PET scan data provides a three-dimensional representation of brain activity with dimensions \(160 \times 160 \times 96\). For this analysis, we select the middle slice along the \(z\)-direction, specifically the 48th slice, to create a two-dimensional domain with \(160 \times 160\) data points, as shown in Fig. \ref{mesh}. We then construct the initial condition \(u_0(x,y)\) based on the PET imaging data. Specifically, we define \(u_0(x,y)\) as
\begin{equation}
u_0(x,y) = \sum_{i=1}^{N_x} \sum_{j=1}^{N_y} u(x_i,y_j) \phi^{(i,j)}_h(x,y), \label{IC}
\end{equation}
where \(\phi^{(i,j)}_h(x,y)\) is the hat function associated with the point \((x_i,y_j)\). As an example, the initial condition for one CN patient is shown in Fig. \ref{mesh}.

We set \( T = 42 \) and initialized \( u_0 \) using the PET images with \( \alpha = 10^6 \), \( \rho = 0.012 \), and \( D = 0.002 \) for Eq.~\eqref{eq:model} and Eq.~\eqref{eq:objective}. We use the linear combination adjoint method and choose \( \beta = 0.5 \) in Algorithm \ref{alg:linear_combination}.

\begin{figure}[htbp]
    \centering
    \begin{tikzpicture}
        \node[anchor=south west, inner sep=0] (image) at (0,0) {\includegraphics[width=0.9\linewidth]{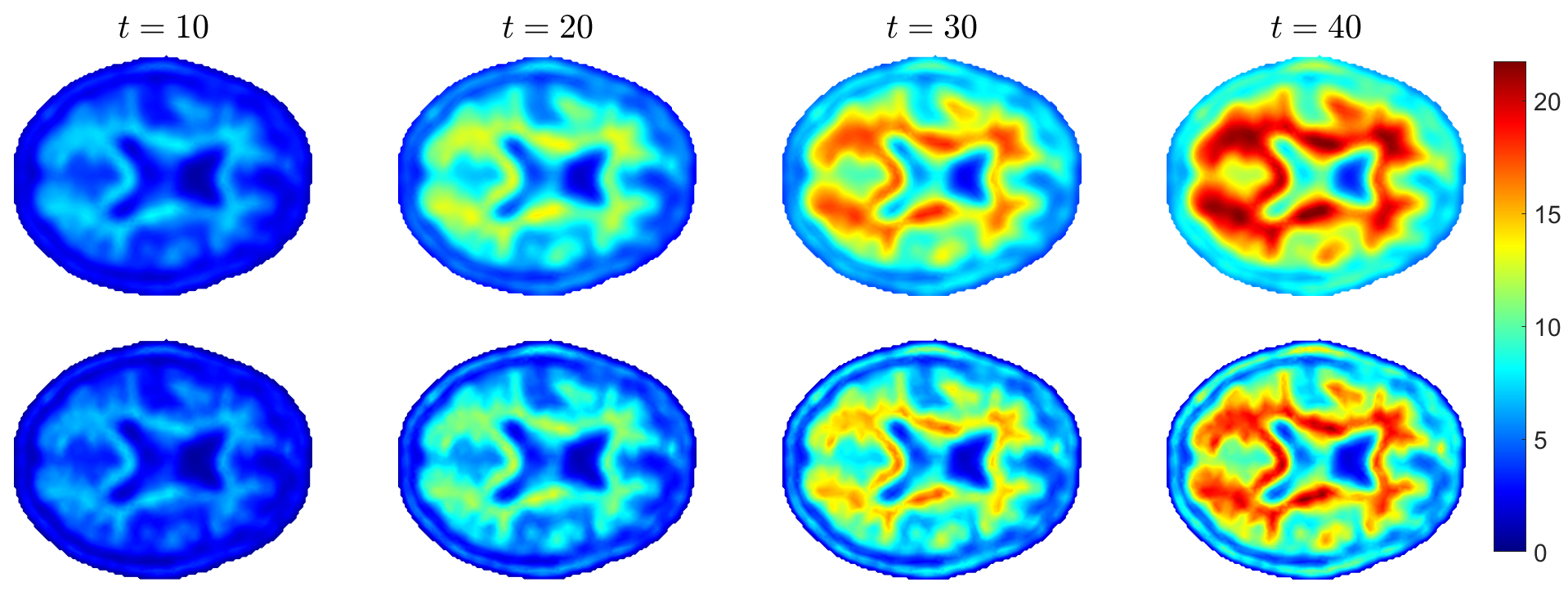}};
        \node[rotate=90] at (-1.2,4.4) {$\int_0^t u(\mathbf{x},s)ds$};
        \node[rotate=90] at (-1.2,1.3) {$\int_0^t u^*(\mathbf{x},s)ds$};
    \end{tikzpicture}
    \caption{The top row shows the temporal accumulation of the constant treatment, represented by \( \int_0^t u(\mathbf{x},s)\,ds \), while the bottom row illustrates the accumulation under the optimal treatment, \( \int_0^t u^*(\mathbf{x},s)\,ds \). Both quantities are displayed at selected time points (\( t = 10, 20, 30, 40 \)), allowing visual comparison of the cumulative effects of the constant and optimal treatments over time.}
    \label{Brain_Figure_U1}
\end{figure}

\begin{figure}
    \centering
    \includegraphics[width=1\linewidth]{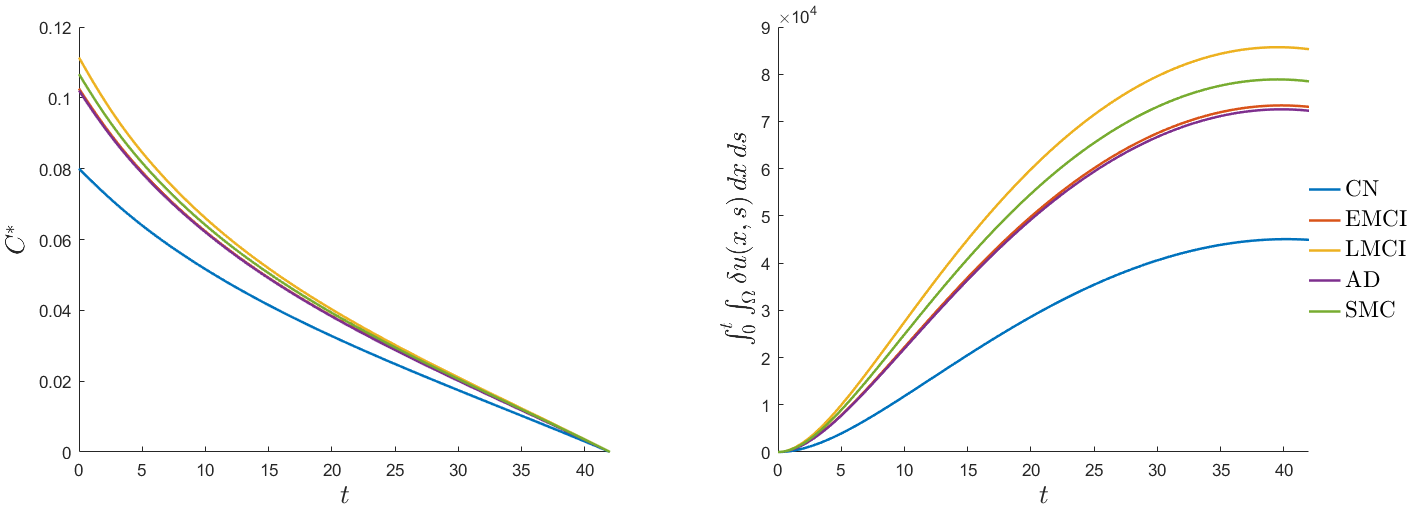}
\caption{{The difference defined as 
$\delta u(\mathbf{x}, t) = u(\mathbf{x}, t) - u^*(\mathbf{x}, t)$, 
where $u$ denotes the solution under the constant treatment $C$, 
and $u^*$ under the optimal treatment $C^*$. For five subjects across five different patient groups: Cognitively Normal (CN), Significant Memory Complaint (SMC) but clinically normal, Early Mild Cognitive Impairment (EMCI), Late Mild Cognitive Impairment (LMCI), and Alzheimer's Disease (AD). We set \( T = 42 \) and initialized \( u_0 \) using the PET images with \( \alpha = 10^6 \), \( \rho = 0.012 \), and \( D = 0.002 \) for Eq.~\eqref{eq:model} and Eq.~\eqref{eq:objective}.}}\label{5groups}
\end{figure}

Similarly to the 1D case, we compare the constant and optimal controls in Fig.~\ref{Brain_Figure_U1}, where \(u\) denotes the constant treatment (Eq.~\eqref{c_avf}) and \(u^*\) the optimal treatment. Figure~\ref{Brain_Figure_U1} presents the temporal evolution of the cumulative integrals \(\int_0^t  u(\mathbf{x},s) \, ds,\) (top row) and \(\int_0^t  u^*(\mathbf{x},s) \, ds,\) (bottom row). The cumulative value under the optimal control remains consistently smaller than that of the constant treatment, indicating that the optimal strategy achieves comparable therapeutic outcomes. By accounting for potential side effects through the cumulative measure \(\int_0^T \int_{\Omega} u(\mathbf{x},s) \, d\mathbf{x} \, ds\), the optimal treatment achieves a better overall balance between amyloid clearance and safety. Figure~\ref{5groups} shows comparisons across five subjects from five patient groups, illustrating that the optimal treatment consistently outperforms the constant treatment under the same parameter settings.

\subsection{Experiments with Empirical Data}

In this section, we perform experiments using empirical abeta PET scan data of the brain surface. The longitudinal PET data allow us to formulate and solve an inverse problem to simultaneously estimate the parameters \(D\) and \(\rho\). Since the brain surface is relatively uniform, we approximate \(D\) as a constant, which suffices to capture the spatiotemporal dynamics of amyloid deposition. Using these estimated parameters, we then apply the optimal control procedure to evaluate treatment strategies under realistic conditions.

We employ PET scans from the Alzheimer's Disease Neuroimaging Initiative (ADNI) dataset \cite{mueller2005alzheimer} (\url{https://adni.loni.usc.edu}). The computational model is built on the brain surface nodes, using a high-resolution mesh consisting of 163,842 nodes and 327,680 triangular elements.

\begin{figure}[htbp]
    \centering
    \includegraphics[width=1\linewidth]{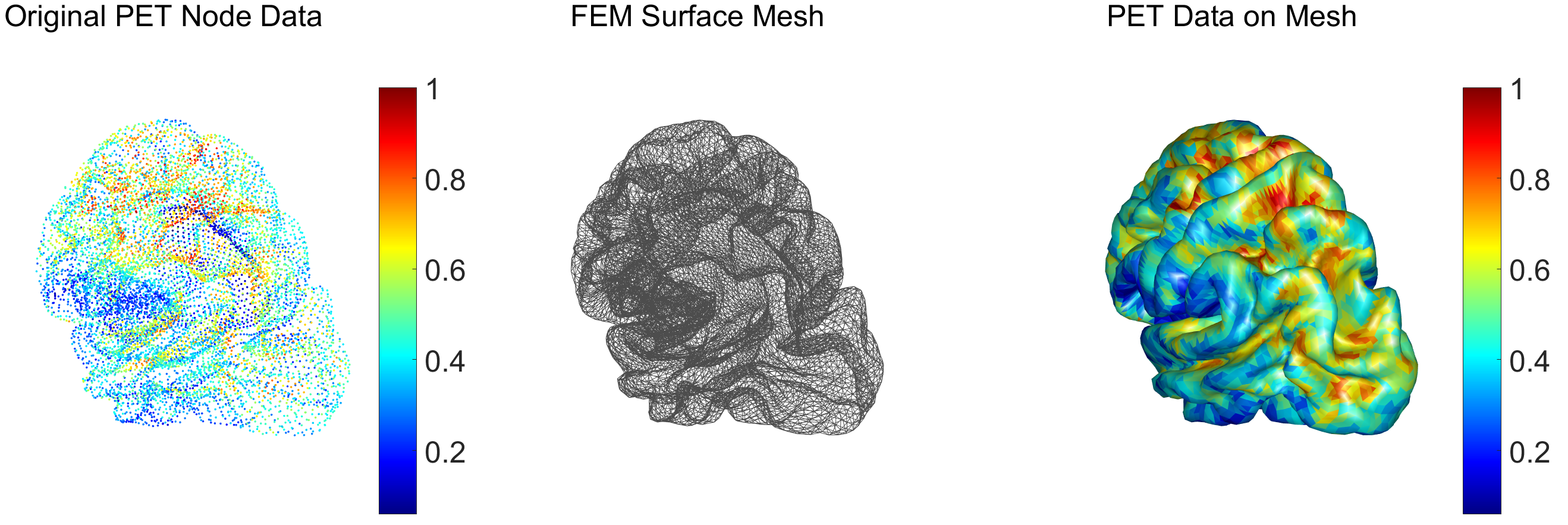}
    \caption{Visualization of the PET-based initial condition. {\bf Left}: Original PET node data visualized as colored scatter points. {\bf Middle}: Surface mesh generated from the brain geometry. {\bf Right}: Initial condition \(u_0\) interpolated over the mesh elements from the PET values.}
    \label{mesh_3D}
\end{figure}

{\bf The initial condition} \(u_0\) is constructed by assigning the baseline PET values to the corresponding surface mesh nodes and representing them using finite elements:
\begin{equation}
u_0(x,y,z) = \sum_{i=1}^{N} u_0(\mathbf{x}_i) \, \phi_i(x,y,z),
\end{equation}
where \(N\) is the number of surface nodes, \(u_0(\mathbf{x}_i)\) is the PET value at node \(i\), and \(\phi_i(\mathbf{x})\) is the associated hat function.
Fig.~\ref{mesh_3D} shows the initial condition construction using a reduced mesh with 10,242 nodes and 20,480 triangular elements. The actual experiments, however, are performed on the high-resolution mesh with 163,842 nodes and 327,680 elements.

\begin{figure}[htbp]
    \centering
    \includegraphics[width=1\linewidth]{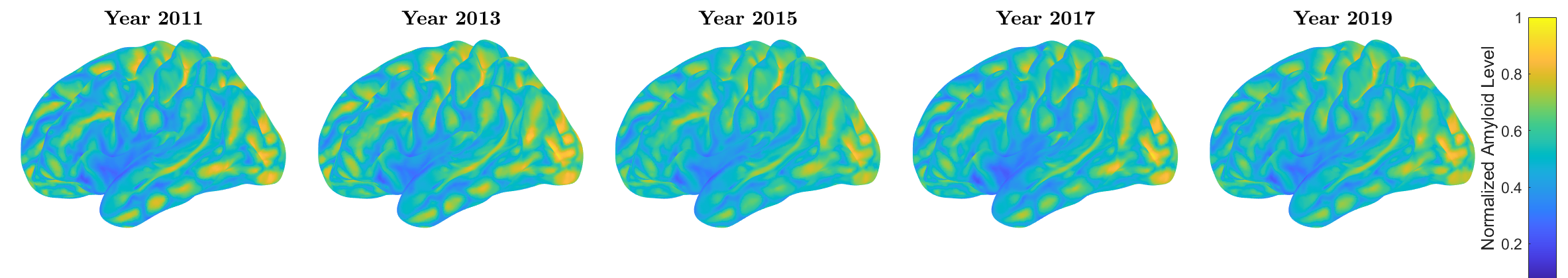}
    \caption{Temporal evolution of normalized amyloid deposition obtained from PET scans, used as observational data for parameter inference.}
    \label{fig:sphere}
\end{figure}

{\bf Parameter Estimation with \texorpdfstring{$\beta$}{Beta}-Amyloid PET Imaging Data}  
We use longitudinal amyloid PET scan data as in Fig.~\ref{fig:sphere} to estimate the parameters \(\rho\) and \(D\) in Eq.~\eqref{eq:model} by minimizing the following training loss:  

\begin{equation}\label{train_loss}
\min_{\rho,D} \frac{1}{N-1}\sum_{i=1}^{N-1} \frac{1}{|\Omega|} \int_{\Omega} 
\frac{|u_{\mathrm{obs}}(\mathbf{x},t_i) - u_{\mathrm{pre}}(\mathbf{x},t_i)|}{u_{\mathrm{obs}}(\mathbf{x},t_i) + \varepsilon} \, d\mathbf{x},
\end{equation}  

where \(u_{\mathrm{obs}}\) is the observed amyloid distribution and \(u_{\mathrm{pre}}\) is the model prediction. The small constant \(\varepsilon = 10^{-4}\) in the denominator prevents numerical instability when \(u_{\mathrm{obs}}\) is close to zero.

Each subject had at least five PET time points. The first \(N-1\) time points were used for training, while the final time point was reserved for validation. The validation loss is defined as  

\begin{equation}\label{vali_loss}
\text{Validation loss} = \frac{1}{|\Omega|} \int_{\Omega} 
\frac{|u_{\mathrm{obs}}(\mathbf{x},t_N) - u_{\mathrm{pre}}(\mathbf{x},t_N)|}{u_{\mathrm{obs}}(\mathbf{x},t_N) + \varepsilon} \, d\mathbf{x}.
\end{equation}

\begin{figure}[htbp]
    \centering
    \includegraphics[width=0.48\linewidth]{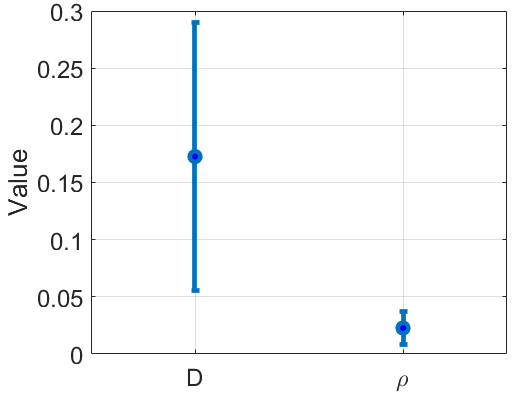}
        \includegraphics[width=0.48\linewidth]{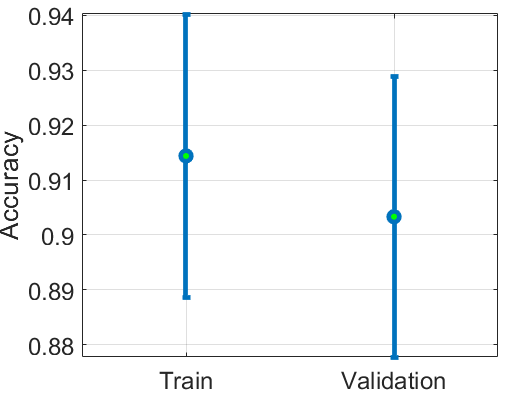}
    \caption{Mean and standard deviation of the diffusion coefficient $D$ and reaction rate $\rho$, and corresponding training and validation accuracy computed as $1 - \text{relative error loss}$ from Eqs.~\eqref{train_loss} and~\eqref{vali_loss}.}
    \label{fig:values}
\end{figure}

The inferred parameter values are summarized in Fig.~\ref{fig:values}, which reports the mean and standard deviation across all subjects. The corresponding training and validation accuracies exceed 90\%, demonstrating the robustness of the model fitting to real PET data.


\begin{figure}[htbp]
    \centering
    \begin{tikzpicture}
        \node[anchor=south west, inner sep=0] (image) at (0,0) {\includegraphics[width=0.9\linewidth]{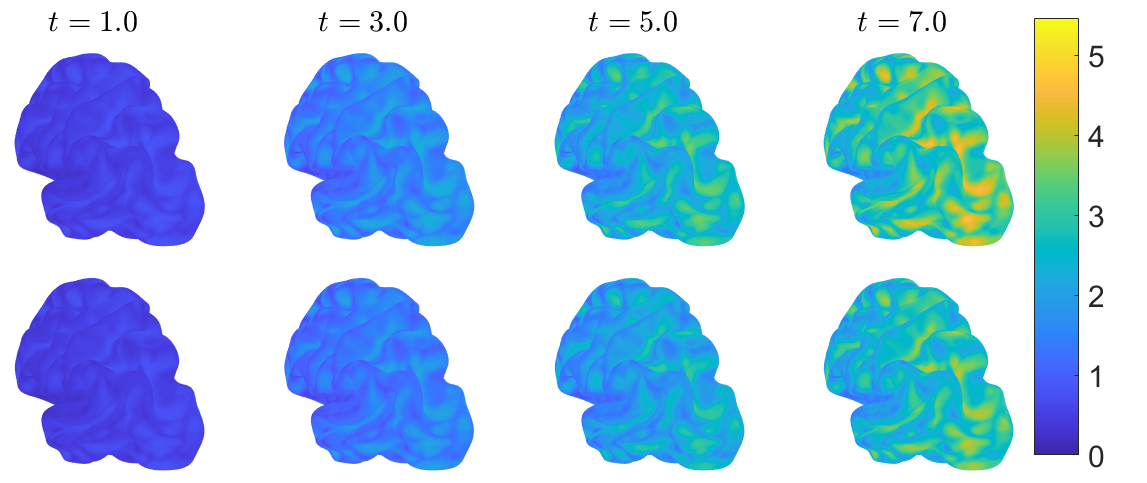}};
        \node[rotate=90] at (-1.2,4.4) {$\int_0^t u(\mathbf{x},s)ds$};
        \node[rotate=90] at (-1.2,1.3) {$\int_0^t u^*(\mathbf{x},s)ds$};
    \end{tikzpicture}
    \caption{The top row shows the temporal accumulation of the constant treatment, represented by \( \int_0^t u(\mathbf{x},s)\,ds \), while the bottom row illustrates the accumulation under the optimal treatment, \( \int_0^t u^*(\mathbf{x},s)\,ds \). Both quantities are displayed at selected time points (\( t = 1, 3, 5, 7 \)), allowing visual comparison of the cumulative effects of the constant and optimal treatments over time.}
    \label{Brain_Figure_3d}
\end{figure}

Once these parameters were estimated, they were incorporated into the optimal control framework to design and evaluate treatment strategies under realistic conditions. This approach enables optimization based on empirically informed model dynamics.
Using the inferred parameters $D$ and $\rho$, we conducted the optimal control experiment over a time horizon of $T = 8$, corresponding to the typical duration of AD progression (eight years). Fig.~\ref{Brain_Figure_3d} illustrates the optimal control results for a representative subject, where the linear combination adjoint method was employed with $\beta = 0.5$ in Algorithm~\ref{alg:linear_combination}. 

\begin{figure}[htbp]
    \centering
    \includegraphics[width=1\linewidth]{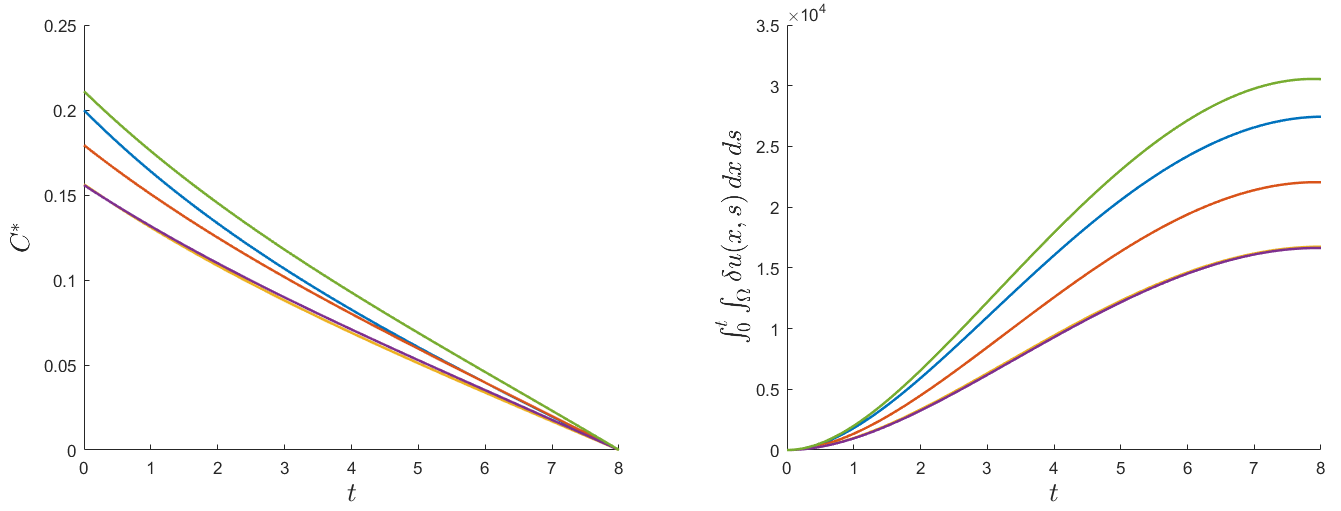}
\caption{{The difference defined as 
$\delta u(\mathbf{x}, t) = u(\mathbf{x}, t) - u^*(\mathbf{x}, t)$, 
where $u$ denotes the solution under the constant treatment $C$, 
and $u^*$ under the optimal treatment $C^*$.} We set $\alpha=5.0e+05$ and $T=8$ for Eq.~\eqref{eq:model} and Eq.~\eqref{eq:objective}. For each subject, $D$ and $\rho$ were individually estimated from their available longitudinal PET scan data, as shown in Fig.~\ref{fig:values}.}\label{all_p}
\end{figure}

We compare the constant control $u$ with the optimal control $u^*$ in Fig.~\ref{Brain_Figure_3d}. Here, $u$ represents a constant treatment applied uniformly over time, whereas $u^*$ denotes the time-dependent optimal control obtained from the optimization procedure. The cumulative difference \(\int_0^t \int_\Omega (u - u^*)\,dx\,ds\) increases during the early and intermediate stages and remains strictly positive throughout the treatment period, indicating a sustained cumulative advantage of the optimal strategy. The overall cumulative effect stays higher, demonstrating the long-term advantage of the optimal control. Unlike the constant treatment, which neither accounts for potential side effects nor provides a clear criterion for determining an appropriate dosage level, the optimal control explicitly incorporates these considerations. By incorporating these effects, the optimal treatment achieves a superior overall balance between amyloid plaque reduction and treatment safety. The results across all subjects are summarized in Fig.~\ref{all_p}. For each subject, the diffusion coefficient $D$ and reaction rate $\rho$ were individually estimated from their available longitudinal PET scan data, and the optimal control procedure was applied using these inferred values, as shown in Fig.~\ref{fig:values}.

\section{Conclusion}

In this study, we have developed a spatially explicit reaction-diffusion framework to optimize anti-amyloid beta therapies for Alzheimer's disease. By incorporating spatial dynamics into a Fisher-KPP-type model and formulating an optimal control problem, we have addressed the critical need to account for the regional heterogeneity of amyloid plaque accumulation in the brain. Our model not only minimizes the A$\beta$ plaque burden but also incorporates a penalty term to mitigate the adverse effects associated with higher treatment doses, such as ARIA, ensuring a balance between therapeutic efficacy and patient safety.

We established a rigorous mathematical foundation for our approach by proving the well-posedness and uniqueness of the optimal control problem. Using the Finite Element Method, we implemented a numerical algorithm—the Linear Combination Adjoint Method—to compute optimal treatment strategies. A significant advancement of this work is the calibration of our model using longitudinal, patient-specific A$\beta$ PET scan data from the ADNI database. This allowed us to infer key biological parameters, such as the plaque growth rate and effective diffusivity, creating a more biologically grounded and personalized model.

Our numerical experiments, conducted in 1D, 2D, and on 3D brain surface geometries, demonstrated the consistent superiority of the optimized treatment strategies over constant dosing regimens. The results highlight that optimized dosing achieves a greater reduction in the cumulative amyloid burden over time, which is crucial for long-term disease management. The application to multiple diagnostic groups—from Cognitively Normal to Alzheimer's Disease—further validated the robustness and generalizability of our approach for personalizing therapeutic interventions.

In conclusion, this study underscores the potential of integrating mathematical modeling, optimal control theory, and clinical imaging data to refine therapeutic strategies for neurodegenerative diseases. By bridging the gap between theoretical modeling and clinical application, our personalized, data-driven framework provides a powerful tool for designing safer and more effective treatment protocols. Future work will extend this framework to incorporate additional biomarkers, such as tau protein, and explore multi-objective optimization to further enhance treatment personalization for Alzheimer's disease.

\section*{Ethics}
This work did not require ethical approval from a human subject or animal welfare committee.

\section*{Data accessibility}
Data used in preparation of this article were obtained from the Alzheimer’s Disease Neuroimaging Initiative (ADNI) database (\url{https://adni.loni.usc.edu}). As such, the investigators within the ADNI contributed to the design and implementation of ADNI and/or provided data but did not participate in analysis or writing of this report. A complete listing of ADNI investigators can be found at: \url{http://adni.loni.usc.edu/wp-content/uploads/how_to_apply/ADNI_Acknowledgement_List.pdf}. The custom MATLAB code developed for the finite element simulations and optimal control analysis in this study is publicly available in the GitHub repository: \texttt{2lineok/Optimal-Control} (\url{https://github.com/2lineok/Optimal-Control}).
\section*{Declaration of AI use}
AI technology has been used to improve the grammar and readability of the text.

\section*{Conflict of interest declaration}
We declare we have no competing interests.
\section*{Acknowlegement}
S.L. and W.H. are supported by both NIH via 1R35GM146894 and NSF DMS-2052685. C. K. acknowledges partial support from NSF DMS 2208373 and DMS 2513176.

\bibliography{new_main_fin.bib}
\appendix
\section{Proofs of the Lemmas}

\subsection{Parabolic Sobolev Embedding}
\begin{lemma}[Parabolic Sobolev Embedding via Interpolation]
\label{thm:parabolic-embedding}
Let $\Omega \subset \mathbb{R}^d$ be a smooth bounded domain, and let 
$Q = \Omega \times (0,T)$ for some $T > 0.$ 
Suppose 
\[
  u \;\in\; L^\infty\bigl(0,T;L^2(\Omega)\bigr) 
  \;\cap\; 
  L^2\bigl(0,T;H^1(\Omega)\bigr).
\]
Then $u$ lies in $L^r(Q)$ for $r = \tfrac{2(d+2)}{d}$, and
\[
  \|u\|_{L^r(Q)} 
  \;\le\; 
  C \Bigl( \|u\|_{L^\infty(0,T;L^2(\Omega))} 
           \;+\; 
           \|\nabla u\|_{L^2(0,T;L^2(\Omega))} \Bigr),
\]
where the constant $C$ depends only on $d,\,T,$ and the domain $\Omega$.
\end{lemma}

\begin{proof}

For a.e.\ $t \in (0,T)$, the function $x \mapsto u(x,t)$ is in $H^1(\Omega)$.  
By the Sobolev embedding \cite{bergh2012interpolation,evans2022partial,triebel1995interpolation}, there is an exponent 
$p \in [2,\,2^*]$ (where $2^* = \tfrac{2d}{\,d-2\,}$ if $d \ge 3$) and a constant $\alpha \in [0,1]$ such that
\[
  \|u(t)\|_{L^p(\Omega)} 
  \;\le\; 
  K_1 \,\|\nabla u(t)\|_{L^2(\Omega)}^{\,\alpha} 
        \,\|u(t)\|_{L^2(\Omega)}^{\,1-\alpha},
\]
where $K_1$ depends only on $\Omega,\,p,$ and $d$, and $\frac{1}{p} 
  \;=\;\alpha\,\Bigl(\tfrac{1}{2}-\tfrac{1}{d}\Bigr) \;+\;(1-\alpha)\,\tfrac{1}{2}.$

We want to prove
$u$ in $L^r(0,T;L^r(\Omega))$.  Take the above estimate, raise to the power $r$, then integrate over $t\in(0,T)$:
\[
  \int_0^T \|u(t)\|_{L^p(\Omega)}^r\,dt
  \;\;\le\;\;
  K_1^r 
  \int_0^T 
     \Bigl[\|\nabla u(t)\|_{L^2(\Omega)}^\alpha 
           \,\|u(t)\|_{L^2(\Omega)}^{\,1-\alpha}
     \Bigr]^r 
  \,dt.
\]
Since $u \in L^\infty(0,T;L^2(\Omega))$, the term $\|u(t)\|_{L^2(\Omega)}^{(1-\alpha)r}$ is bounded in $t$, so we can factor it out of the integral. 
It remains to check the integrability of $\|\nabla u(t)\|_{L^2(\Omega)}^{\,\alpha r}$.  
Because $\nabla u \in L^2(0,T;L^2(\Omega))$, we need
$\alpha\,r \;\le\;2$ to ensure that $\bigl\|\nabla u(t)\bigr\|_{L^2(\Omega)}^{\,\alpha r}$ is integrable over $(0,T)$.  

To prove $u \in L^r(Q)$ with the same exponent $r$ in both space and time, we set $p=r.$  
From the Gagliardo--Nirenberg relation \cite{evans2022partial} we have $\frac{1}{p} 
  \;=\;\frac{1}{2} - \frac{\alpha}{d}.$
Hence if $p = r$, then 
\[
  \frac{1}{r} 
  \;=\;
  \frac{1}{2} - \frac{\alpha}{d} 
  \;\;\Longrightarrow\;\; 
  r \;=\;\frac{1}{\tfrac12 - \tfrac{\alpha}{d}} 
        \;=\;\frac{2d}{d - 2\alpha}.
\]
The time‐integrability requirement $\alpha\,r \le 2$ becomes
\[
  \alpha \, \frac{2d}{d-2\alpha} 
  \;\;\le\; 
  2
  \;\;\Longrightarrow\;\;
  \alpha(d+2) \;=\; d
  \;\;\Longrightarrow\;\; 
  \alpha = \frac{d}{d+2}.
\]
Substitute $\alpha = \tfrac{d}{d+2}$ back to find $ r 
  \;=\; 
  \frac{2(d+2)}{d}.$
Hence $u \in L^r(Q)$ with $r = \tfrac{2(d+2)}{d}$.  Moreover, with young's inequality we have
\[
  \|u\|_{L^r(Q)} 
  \;\le\;
  C \Bigl(\|u\|_{L^\infty(0,T;L^2(\Omega))}
        + \|\nabla u\|_{L^2(0,T;L^2(\Omega))}\Bigr),
\]
where $C$ depends only on $\Omega,\,d,$ and $T$. 
This completes the proof.
\end{proof}

\subsection{Derivative of the Equation}
\begin{lemma} \label{thm:necessary_condition_proof}
    Let \( C \in L^{\infty}(0,T) \) be non-negative, and let \( u = u_C \) be the corresponding solution of \eqref{eq:model}. For \( \eta \in L^{\infty}(0,T) \), the mapping \( C \to u_C(\mathbf{x}, t) \) is differentiable in the following sense: there exists a function \( \psi = \psi_{C, \eta} \in L^2(0,T; H^1(\Omega)) \) such that
    \begin{equation*}
        \psi_\epsilon \rightharpoonup \psi \quad \text{weakly in } L^2(0, T; H^1(\Omega)) \quad \text{as } \epsilon \to 0,
    \end{equation*}
    where \( \psi_\epsilon = \psi_{C, \epsilon, \eta} = \frac{u_\epsilon - u}{\epsilon} \), \( u_\epsilon = u_{C + \epsilon \eta} \), and the sensitivity \( \psi \) is the weak solution of
    \begin{equation} \label{eq:sensitivity_proof}
        \begin{cases}
            \psi_{t} - \nabla \cdot \big( D(\mathbf{x}) \nabla \psi \big) - (\rho - 2 \rho u - C) \psi = -\eta u, & \text{in } \Omega \times (0, T), \\[1ex]
            \frac{\partial \psi}{\partial \mathbf{n}} = 0, & \text{on } \partial \Omega \times (0, T), \\[1ex]
            \psi(\mathbf{x}, 0) = 0, & \text{in } \Omega.
        \end{cases}
    \end{equation}

    Moreover, the objective functional satisfies the expansion:
    \begin{equation} \label{eq:obj_expansion_proof}
        \mathcal{J}(C + \epsilon \eta) = \mathcal{J}(C) + \epsilon \int_0^T \left( \int_{\Omega} \psi \, d\mathbf{x} + 2 \alpha \eta C \right) dt + o(\epsilon).
    \end{equation}
\end{lemma}

\begin{proof}
\noindent
We follow the idea of the proof to \cite{yousefnezhad2021optimal}. In the following proof, \( K\) denotes a generic positive constant independent of \( \epsilon \). By Eq.~\eqref{eq:model} and Lemma \ref{thm:existence}, \( u_{\epsilon} \) is the unique weak solution of the following equation:
\begin{equation} \label{eq:uepsilon}
    \begin{cases}
        \partial_t u_{\epsilon} - \nabla \cdot \left( D(\mathbf{x}) \nabla u_{\epsilon} \right) = \rho (1 - u_{\epsilon}) u_{\epsilon} - (C + \epsilon \eta) u_{\epsilon}, & \text{in } \Omega \times (0, T), \\
        \frac{\partial u_{\epsilon}}{\partial \mathbf{n}} = 0, & \text{on } \partial \Omega \times (0, T), \\
        u_{\epsilon}(\mathbf{x}, 0) = u_0(\mathbf{x}), & \text{in } \Omega.
    \end{cases}
\end{equation}
From Lemma \ref{thm:uniform}, we know that
\begin{equation*}
    \|u_{\epsilon}\|_{L^{\infty}(Q_T)} + \|\partial_t u_{\epsilon}\|_{L^2(0, T; H^1(\Omega)^*)} + \|u_{\epsilon}\|_{L^2(0, T; H^1(\Omega))} \leq K,
\end{equation*}
where \( K \) is a positive constant for all \( 0<\epsilon<1 \).

Define \( \psi_{\epsilon} = \frac{u_{\epsilon} - u}{\epsilon} \), where \( u \) solves \eqref{eq:model} with control \( C \). Subtracting \eqref{eq:model} from \eqref{eq:uepsilon}, we find that \( \psi_{\epsilon} \) solves
\begin{equation} \label{eq:psi_epsilon}
    \begin{cases}
        \partial_t \psi_{\epsilon} - \nabla \cdot \left( D(\mathbf{x}) \nabla \psi_{\epsilon} \right) - (\rho - \rho u_{\epsilon} - \rho u - C) \psi_{\epsilon} = -\eta u_{\epsilon}, & \text{in } \Omega \times (0, T), \\
        \frac{\partial \psi_{\epsilon}}{\partial \mathbf{n}} = 0, & \text{on } \partial \Omega \times (0, T), \\
        \psi_{\epsilon}(\mathbf{x}, 0) = 0, & \text{in } \Omega.
    \end{cases}
\end{equation}
By standard results of linear parabolic PDEs (e.g., Theorem 1.1.2 in \cite{reichelt2015two}), there exists a unique weak solution to \eqref{eq:psi_epsilon}.

We now show that as \( \epsilon \to 0 \), \( \psi_{\epsilon} \rightharpoonup \psi \) weakly in \( L^2(0, T; H^1(\Omega)) \), where \( \psi \) solves \eqref{eq:sensitivity}.

Multiplying both sides of \eqref{eq:psi_epsilon} by \( \psi_{\epsilon} \) and integrating over \( \Omega \), we obtain:
\begin{equation*}
    \frac{1}{2} \frac{d}{dt} \int_{\Omega} \psi_{\epsilon}^2 \, d\mathbf{x} + \int_{\Omega} D(\mathbf{x}) |\nabla \psi_{\epsilon}|^2 \, d\mathbf{x} 
    = \int_{\Omega} (\rho - \rho u_{\epsilon} - \rho u - C) \psi_{\epsilon}^2 \, d\mathbf{x} - \int_{\Omega} \eta u_{\epsilon} \psi_{\epsilon} \, d\mathbf{x}.
\end{equation*}
Integrating over \( (0, t) \), we get:
\small\begin{equation} \label{eq:energy_estimate}
    \frac{1}{2} \|\psi_{\epsilon}(t)\|_{L^2(\Omega)}^2 + \int_0^t \int_{\Omega} D(\mathbf{x}) |\nabla \psi_{\epsilon}|^2 \, d\mathbf{x} ds \leq \rho \int_0^t \|\psi_{\epsilon}\|_{L^2(\Omega)}^2 \, ds + \|u_{\epsilon}\|_{L^{\infty}(Q_T)} \int_0^t \int_{\Omega} |\eta \psi_{\epsilon}| \, d\mathbf{x} ds.
\end{equation}
Using the Cauchy-Schwarz inequality and the bound \( 2 |\eta \psi_{\epsilon}| \leq \eta^2 + \psi_{\epsilon}^2 \), we find:
\begin{equation*}
    \|\psi_{\epsilon}(t)\|_{L^2(\Omega)}^2 \leq \left( 2\rho + \|u_{\epsilon}\|_{L^{\infty}(Q_T)} \right) \int_0^t \|\psi_{\epsilon}\|_{L^2(\Omega)}^2 \, ds + \|u_{\epsilon}\|_{L^{\infty}(Q_T)} |\Omega| \|\eta\|_{L^2(0, T)}.
\end{equation*}
By Grönwall's inequality \cite{evans2022partial}, we conclude that \( \|\psi_{\epsilon}\|_{L^\infty(0, T; L^2(\Omega))} \leq K \). Combining this with the bound on the gradient term in \eqref{eq:energy_estimate}, we get:
\begin{equation} \label{eq:psi_H1_bound}
    \|\psi_{\epsilon}\|_{L^2(0, T; H^1(\Omega))} \leq K.
\end{equation}

Multiplying both sides of \eqref{eq:psi_epsilon} by \( \phi \in L^2(0, T; H^1(\Omega)) \) and integrating over \( Q_T \), we obtain:
\begin{equation*}
    \begin{aligned}
        \int_0^T \int_{\Omega} (\partial_t \psi_{\epsilon}) \phi \, d\mathbf{x} dt &= - \int_0^T \int_{\Omega} D(\mathbf{x}) \nabla \psi_{\epsilon} \cdot \nabla \phi \, d\mathbf{x} dt \\ &\
        - \int_0^T \int_{\Omega} \eta u_{\epsilon} \psi_{\epsilon} \phi \, d\mathbf{x} dt + \int_0^T \int_{\Omega} (\rho - \rho u_{\epsilon} - \rho u - C) \psi_{\epsilon} \phi \, d\mathbf{x} dt.
    \end{aligned}
\end{equation*}
Using the bound-ness for \( \psi_{\epsilon} \), we conclude:
\begin{equation} \label{eq:psi_dt_bound}
    \|\partial_t \psi_{\epsilon}\|_{L^2(0, T; H^1(\Omega)^*)} \leq K.
\end{equation}

By passing to a subsequence, it follows that
\begin{equation*}
    \begin{aligned}
    & u_{\epsilon} \rightharpoonup u \text{ in } L^2\left(0, T ; H^1(\Omega)\right)\\
        & u_{\epsilon} \rightarrow u \text{ in } L^2(Q_T)\\
        &\partial_t u_{\epsilon} \rightharpoonup \partial_t u \text{ in } L^2(0, T ; H^1(\Omega)^*)\\
        & \psi_{\epsilon} \rightharpoonup \psi \text{ in } L^2\left(0, T ; H^1(\Omega)\right)\\
        &\partial_t \psi_{\epsilon} \rightharpoonup \partial_t \psi \text{ in } L^2(0, T ; H^1(\Omega)^*)
    \end{aligned}
\end{equation*}
Similar to the proof in Thm 3.1 in \cite{yousefnezhad2021optimal} and use the uniqueness of the weak solution $u_C$,  we can conclude that $\psi$ is the weak solution of \eqref{eq:sensitivity}
\end{proof}

\subsection{Adjoint Equation}
\begin{lemma} \label{thm:adjoint_solution_proof}
    Let \( C \in L^{\infty}(0,T) \) be non-negative, and let \( u = u_C \) be the corresponding solution of \eqref{eq:model}. Then, there exists \( w = w_C \in L^2(0,T; H^1(\Omega)) \) with \( w_t \in L^2(0,T; H^1(\Omega)^*) \) such that \( w \) is the weak solution of
    \begin{equation} \label{eq:adjoint_equation_proof}
        \begin{cases}
            w_{t} + \nabla \cdot \big( D(\mathbf{x}) \nabla w \big) + (\rho - 2 \rho u - C) w = 1, & \text{in } \Omega \times (0, T), \\
            \frac{\partial w}{\partial \mathbf{n}} = 0, & \text{on } \partial \Omega \times (0, T), \\
            w(\mathbf{x}, T) = 0, & \text{in } \Omega.
        \end{cases}
    \end{equation}
    Moreover, there exists a constant \( K > 0 \), depending on depending only on \( |\Omega| \), \( T \), \( \rho \), \( \theta \), \( \|u_0\|_{L^{\infty}(\Omega)} \), \( \|C\|_{L^{\infty}(0,T)} \), and spatial dimension $d$, such that:
    \begin{equation*}
        \|w\|_{L^{\infty}(0, T; L^2(\Omega))} + \|w_t\|_{L^2(0, T; H^1(\Omega)^*)} + \|w\|_{L^2(0, T; H^1(\Omega))} \leq K.
    \end{equation*}
    Additionally, we can confirm \( w \leq 0 \) almost everywhere, and there exists a positive constant \( M \), depending only on \( |\Omega| \), \( T \), \( \rho \), \( \theta \), \( \|u_0\|_{L^{\infty}(\Omega)} \), \( \|C\|_{L^{\infty}(0,T)} \), and spatial dimension $d$, such that
    \begin{equation*}
        \|w\|_{L^{\infty}(Q_T)} \leq M.
    \end{equation*}
    
\end{lemma}

\begin{proof}

We follow the idea of the proof to \cite{finotti2012optimal} and \cite{yousefnezhad2021optimal}. Define \( v(\mathbf{x}, t) := -w(\mathbf{x}, T - t) \). Then \( v \) is the unique weak solution of the following problem by standard results of linear parabolic PDEs (Theorem 1.1.2 in \cite{reichelt2015two}):
\begin{equation} \label{eq:dual_equation}
    \begin{cases}
        v_{t} - \nabla \cdot \big( D(\mathbf{x}) \nabla v \big) - \big( \rho - 2 \rho u(\mathbf{x}, T-t) - C(T-t) \big) v = 1, & \text{in } \Omega \times (0, T), \\
        \frac{\partial v}{\partial \mathbf{n}} = 0, & \text{on } \partial \Omega \times (0, T), \\
        v(\mathbf{x}, 0) = 0, & \text{in } \Omega.
    \end{cases}
\end{equation}

Multiplying both sides of \eqref{eq:dual_equation} by \( v \) and integrating over \( \Omega \), we obtain:
\small\begin{equation*}
    \frac{1}{2} \frac{d}{dt} \int_{\Omega} v^2 \, d\mathbf{x} + \int_{\Omega} D(\mathbf{x}) |\nabla v|^2 \, d\mathbf{x} + \int_{\Omega} \big( 2 \rho u(\mathbf{x}, T-t) + C(T-t) \big) v^2 \, d\mathbf{x} = \int_{\Omega} \rho v^2 \, d\mathbf{x} + \int_{\Omega} v \, d\mathbf{x}.
\end{equation*}
Applying the Cauchy-Schwarz inequality with knowing $D(x)|\nabla v|^2+(2\rho u +C)v^2\geq 0$, we obtain:
\begin{equation*}
    \frac{d}{dt} \int_{\Omega} v^2 \, d\mathbf{x} \leq (2 \rho + 1) \int_{\Omega} v^2 \, d\mathbf{x} + |\Omega|.
\end{equation*}
By Grönwall's inequality \cite{evans2022partial}, we conclude that:
\begin{equation*}
    \|v\|_{L^{\infty}(0, T; L^2(\Omega))} \leq K.
\end{equation*}
The corresponding bounds on \( \nabla v \) and \( v_t \) can be established similarly, completing the proof of first inequality.

Now we will show the $\|\cdot\|_{L^\infty(Q_T)}$ bound. For any fixed $K \in \mathbb{N}$, let $0 = t_0 < t_1 < \cdots < t_K = T$ be a partition of $[0, T]$ which will be determined. For each $i = 1, 2, \cdots, K$ let $Q_i = \Omega \times [t_{i-1}, t_i]$ and
\begin{equation*}
    \|v\|_{V_2(Q_i)}^2 := \sup_{t \in [t_{i-1}, t_i]} \int_\Omega v^2(x,t)dx + \int_{Q_i} |\nabla v(x,t)|^2 dxdt.
\end{equation*}
By observation, we have $v \geq 0$. Therefore, we only need to show that $v$ is bounded above. It suffices to show that $\|v\|_{L^\infty(Q_i)}$ is bounded above for all $i = 1,2,\dots,K$.

By previous computations, we see that there exists a finite constant $C_0$ (depends only  \( |\Omega| \), \( T \), \( \rho \), \( \theta \), \( \|u_0\|_{L^{\infty}(\Omega)} \), and \( \|C\|_{L^{\infty}(0,T)} \)) such that
\begin{equation}
    \|v\|_{V_2(Q_i)} \leq C_0, \quad \forall \; i = 1,2,\dots,K.
\end{equation}

Next, for each $k > \hat{k} =\|v_0\|_{L^\infty(\Omega)} + 1$, let us denote
\begin{equation*}
    v^{(k)}(x,t) := \max\{v(x,t) - k, 0\}.
\end{equation*}

Also, denote the sets
\begin{equation*}
    A_k(t) := \{x \in \Omega : v(x,t) > k\}, \quad Q_i(k) := \{(x,t) \in Q_i : v(x,t) > k\}, \quad i = 1,\dots,K.
\end{equation*}

Multiplying the first equation of (2.1) by $v^{(k)}$ and using integration by parts, we get
\small\begin{equation*}
    \frac{1}{2} \frac{d}{dt} \int_\Omega v^{(k)}(x,t)^2 dx +  \int_\Omega D(x)|\nabla v^{(k)}|^2 dx = \int_\Omega [ v v^{(k)}(\rho - 2\rho u(\mathbf{x}, T-t) - C(T-t))+v^{(k)}]dx.
\end{equation*}
Rewriting,
\footnotesize\begin{equation*}
    \int_\Omega [ v v^{(k)}(\rho - 2\rho u(\mathbf{x}, T-t) - C(T-t))]dx \leq \int_{A_k(t)} [(\rho +2\rho \|u\|_{L^\infty(Q_T)}+\|C\|_{L^\infty(0,T)})vv^{(k)} ]dx.
\end{equation*}
\begin{equation*}
    \int_\Omega 2v^{(k)}dx \leq \int_{A_k(t)} [1+v^{(k)}v^{(k)} ]dx  \leq \int_{A_k(t)} [k^2+v^{(k)}v^{(k)} ]dx.
\end{equation*}
We obtain
\begin{equation*}\label{C2dfine}
      \frac{1}{2} \frac{d}{dt} \int_\Omega v^{(k)}(x,t)^2 dx +  \int_\Omega D(x)|\nabla v^{(k)}|^2 dx \leq C_2 \int_{A_k(t)} [(v - k)^2 + k^2] dx..
\end{equation*}

for some constants $C_2 > 0$ depending only on \( |\Omega| \), \( T \), \( \rho \), \( \theta \), \( \|u_0\|_{L^{\infty}(\Omega)} \), and \( \|C\|_{L^{\infty}(0,T)} \). Note that $u^{(k)}(\cdot,0) = 0$. Thus, by integrating this equation in time on $[0,t]$ with $0 < t < t_1$, we obtain
\begin{equation*}
    \|v^{(k)}\|_{V_2(Q_i)}^2 \leq C_2 \int_{Q_i(k)} [(v - k)^2 + k^2]dxdt.
\end{equation*}

Note that
\small\begin{equation}
\int_{Q_1(k)} (v - k)^2 \,dx\,dt 
\;=\;\int_{Q_1(k)} \bigl[v^{(k)}\bigr]^2 \,dx\,dt
\;\le\; 
t_1 \sup_{0 < t < t_1} \int_{\Omega} \bigl[v^{(k)}(x,t)\bigr]^2 \,dx 
\;\le\;
t_1 \,\bigl\|v^{(k)}\bigr\|_{V_2(Q_1)}^2.
\end{equation}
Therefore, choosing $t_1$ sufficiently small such that $t_1C_2 < \tfrac12$ yields
\begin{equation}
\bigl\|v^{(k)}\bigr\|_{V_2(Q_1)}^2 
\;\le\;
2C_2\,k^2\,\sigma(k),
\quad
\text{where}
\quad
\sigma(k) \;:=\; |Q_1(k)|
\;=\;\int_{0}^{t_1} |A_k(t)| \,dt.
\end{equation}
Equivalently,
\begin{equation}\label{sresult}
\bigl\|v^{(k)}\bigr\|_{V_2(Q_1)}
\;\le\;
C_3\,k\,\bigl[\sigma(k)\bigr]^{\tfrac12},
\quad
\forall\,k > \hat{k}.
\end{equation}

First of all, for all $2 \;\le\; r \;\le\; \frac{2(d+2)}{d},$ by the Sobolev embedding in Lemma \ref{thm:parabolic-embedding}, we can find a constant $\beta_0>0$ (depending only on $|\Omega|,d,r,$ and $T$) such that
\begin{equation}
\|w\|_{L^r(Q_k)} \;\le\;\beta_0\,\|w\|_{V_2(Q_k)},
\quad
\forall\,w\in V_2(Q_k), \;\;\forall\,k=1,2,\dots,K.
\end{equation}

Let $M_0 = m_0\,\hat{k}$ for some $m_0>1$ which will be determined later. Also, for $i=0,1,2,\dots,$ let us denote 
\[
k_i := M_0 (2 - 2^{-i}).
\]
It follows directly from the definition of $\sigma$ that
\begin{equation}
\bigl(k_{i+1} - k_i\bigr)\,\sigma^\frac{1}{r}\bigl(k_{i+1}\bigr)
\;\le\;\bigl\|v^{(k_i)}\bigr\|_{L^r(Q_1)},
\quad
\forall\,i\in\mathbb{N}\cup\{0\}.
\end{equation}
From now on, we fix
\[
2 \;<\; r \;<\; \frac{2(d+2)}{d}
\]
and write $r = 2(1 + \kappa)$ for some $\kappa>0$. Since $k_i>\hat{k}$ for all $i$, from \eqref{sresult} we have
\begin{equation}
\bigl\|v^{(k_i)}\bigr\|_{L^r(Q_1)}
\;\le\;\beta_0\,\bigl\|v^{(k_i)}\bigr\|_{V_2(Q_1)}
\;\le\;\beta_0\,C_3\,k_i\,\bigl[\sigma(k_i)\bigr]^{\tfrac{1+\kappa}{r}},
\quad
\forall\,i\in\mathbb{N}\cup\{0\}.
\end{equation}
Then, combining inequalities, we get
\begin{equation}
\sigma\bigl(k_{i+1}\bigr)^{\tfrac{1}{r}}
\;\le\;\frac{\beta_0\,C_3\,k_i}{\,k_{i+1}-k_i\,}\,
\bigl[\sigma(k_i)\bigr]^{\tfrac{1+\kappa}{r}}
\;\le\;4\,\beta_0\,C_3\,2^i \,\bigl[\sigma(k_i)\bigr]^{\tfrac{1+\kappa}{r}},
\quad
\forall\,i\in\mathbb{N}\cup\{0\}.
\end{equation}

For all $i=0,1,\dots,$ let $y_i = \bigl[\sigma(k_i)\bigr]^{\tfrac{1}{r}}.$ Then it follows directly from the recursion formula and a straightforward induction that
\[
y_i \;\le\; \bigl[4\,\beta_0\,C_3\bigr]^{\frac{(1+\kappa)^i-1}{\kappa}}\,
2^{\frac{(1+\kappa)^i-1}{\kappa^2}-\frac{i}{\kappa}}
\,y_0^{\,(1+\kappa)\,i},
\quad
\forall\,i=0,1,2,\dots.
\]
By similar computation we have
\[
\sigma(M_0)^{\tfrac{1}{r}}
\;\le\;
\frac{\beta_0\,C_3}{\,m_0 - 1\,}
\,\bigl[\sigma(k_i)\bigr]^{\tfrac{1+\kappa}{r}}
\;\le\;
\frac{\beta_0\,C_3}{\,m_0 - 1\,}\bigl[T\,|\Omega|\bigr]^{\tfrac12}.
\]

Thus, by choosing
\[
  m_{0} 
  \;=\; 
  1 
  \;+\; 
  \beta_{0}\,C_3\,\bigl[T\,|\Omega|\bigr]^{\tfrac12}
  \,\bigl(4\,\beta_{0}\,C_3\bigr)^{\tfrac{1}{\kappa}}
  \,2^{\tfrac{1}{\kappa^{2}}},
\]
we have
\[
  y_{0}
  \;=\;
  \sigma(k_{0})^{\tfrac{1}{\kappa}}
  \;=\;
  \sigma(M_{0})^{\tfrac{1}{\kappa}}
  \;\;\le\;\;
  \bigl(4\,\beta_{0}\,C_3\bigr)^{-\tfrac{1}{\kappa}}
  \,2^{-\tfrac{1}{\kappa^{2}}}.
\]

Then, it follows from the inequalities 
\[
  y_{i}
  \;\le\;
  \bigl[\,4\,\beta_{0}\,C_3\bigr]^{-\tfrac{1}{\kappa}}
  \,2^{-\tfrac{1}{\kappa}}
  \,2^{-\tfrac{i}{\kappa}},
  \quad
  \forall\,i = 0,1,2,\dots.
\]
In particular, 
\[
  y_{i} 
  \;=\;
  \sigma\bigl(k_{i}\bigr)^{\tfrac{1}{\kappa}}
  \;\longrightarrow\;
  0
  \quad\text{as } i\to\infty.
\]
Hence, $\sigma\bigl(2M_{0}\bigr)=0$ and therefore, on $Q_{1}$,
\[
  v \;\le\; c_{1}
  \; \stackrel{:=}\;
  2\,m_{0}\,\hat{k}
  \;=\;
  2\,\Bigl\{
        1 
        + 
        \beta_{0}\,C_3\,\bigl[T\,|\Omega|\bigr]^{\tfrac12}
        \,\bigl(4\,\beta_{0}\,C_3\bigr)^{\tfrac{1}{\kappa}}
        \,2^{\tfrac{1}{\kappa^{2}}}
     \Bigr\}
  \,\Bigl\{
       \|v_{0}\|_{L^\infty(Q_{T})} \;+\; 1
     \Bigr\}.
\]

Next, note that similarly to the choice of $t_1$, we choose $K \in \mathbb{N}$ sufficiently large so that
\begin{equation}
C_{2}\,\bigl|t_{k} - t_{k-1}\bigr| 
\;<\;\tfrac{1}{2}
\quad\text{and}\quad
k = 2,\,\dots,\,K,
\end{equation}
where $C_{2}$ is defined in Eq.~\eqref{C2dfine}. Therefore, by the same proof as before, but using $v(\cdot,t_1)$ as $v_{0}$, we can prove that $v$ is bounded above on $Q_{2}$ by some constant $c_{2}$. Repeating this argument iteratively, we arrive at
\[
  \sup_{Q_{i}} v(x,t) \;\le\; c_{i},
  \quad\text{for all}\; i = 2,3,\dots,K,
\]
where all of the constants $c_{i}$ can be explicitly defined as
\begin{equation}
c_{i} 
\;=\; 
2\,\Bigl\{\,1 
  \;+\; \beta_{0}\,C_{3}\,\bigl[T\,|\Omega|\bigr]^{\tfrac12}
         \,\bigl(4\,\beta_{0}\,C_{3}\bigr)^{\tfrac1\kappa}
         \,2^{\tfrac1{\kappa^{2}}}
  \Bigr\}
  \,\bigl(c_{i-1} + 1\bigr).
\tag{3.15}
\end{equation}

Moreover, we see we can choose $K$ large enough so that $K > 2\,T\,C_{2}$. All the constants depends only on $\|C\|_{L^{\infty}(0,T)},\theta,M,|\Omega|,T,\|u_{0}\|_{L^\infty}$, and the spatial dimension $d$. Therefore,
\[
  \sup_{Q_{T}} v \;\le\; C
  \quad\text{with}\quad
  C \;=\; c_{K}
  \;=\;\max\{\,c_{i}\,\mid\,i = 1,2,\dots,K\}.
\]
The proof of the Lemma is therefore complete.
\end{proof}

\end{document}